\documentclass [reqno,  12pt, a4paper]{amsproc}

\usepackage[colorlinks]{hyperref}
\hypersetup{colorlinks,breaklinks,
            urlcolor={blue!0!black},
            linkcolor={black!100!blue},
            citecolor={black!100!blue}
}

\usepackage[margin=1.05in] {geometry}
\usepackage{amsthm,array,mathrsfs,amscd}
\usepackage{amsmath,amssymb,enumitem}
\usepackage{dsfont}

\usepackage[bbgreekl]{mathbbol}

\usepackage{datetime}

\usepackage{tikz}
\usepackage{tikz-cd}
\usetikzlibrary{trees}
\usetikzlibrary{patterns}
\usetikzlibrary{calc}
\usetikzlibrary{decorations.markings}
\usepackage[linguistics]{forest}

\usepackage{booktabs}

\def\half#1#2{\begin{matrix}\frac{#1}{#2}\end{matrix}}

\def\M{\mathscr{M}}

\def\field{\mathbb{k}}
\def\fie{\mathbb{k}}

\def\alg{\mathbb{A}}
\def\balg{\mathbb{B}}

\newcommand{\fusion}{\mathscr{F}}
\newcommand{\basis}{\mathscr{B}}
\newcommand{\mult}{\cdot}

\newcommand{\Id}{\operatorname{Id}}

\newcommand{\nass}{\operatorname{\mathsf{NC}}}
\newcommand{\Nring}[2]{#1 (#2)}

\newcommand{\evan}{\mathcal{E}}

\newcommand{\medmlt}{\odot}

\DeclareMathOperator{\Tri}{\mathcal{S}}

\DeclareMathOperator{\im}{im} 

\newcommand{\chr}{\operatorname{char}}
 
\DeclareMathOperator{\Idm}{Idm}

\DeclareMathOperator{\End}{End}

\newcommand{\NCA}{\operatorname{NCA}}

\DeclareMathOperator{\weight}{deg}
\DeclareMathOperator{\height}{ht}

\def\half#1#2{\begin{matrix}\frac{#1}{#2}\end{matrix}}

\def\M{\mathscr{M}}

\def\len1{1.5cm}
\def\len2{2cm}
\def\XX{$\bullet$}
\def\NN{node {\XX}}

\def\alg{\mathbb{A}}

\newtheorem{theorem}{\rm\bf Theorem}[section]
\newtheorem{proposition}[theorem]{\rm\bf Proposition}
\newtheorem{lemma}[theorem]{\rm\bf Lemma}
\newtheorem{corollary}[theorem]{\rm\bf Corollary}
\theoremstyle{definition}
\newtheorem{definition}[theorem]{\rm\bf Definition}

\theoremstyle{remark}
\newtheorem{remark}[theorem]{\rm\bf Remark}
\newtheorem{example}[theorem]{\rm\bf Example}

\begin{document}

\title[Commutative algebras satisfying]{Commutative algebras satisfying univariate identities with vanishing Peirce polynomial}

\author[D.~J.~F.~Fox]{Daniel J. F. Fox}
\address[D.~J.~F.~Fox]{Departamento de Matemática Aplicada\\ Escuela Técnica Superior de Arquitectura\\ Universidad Politécnica de Madrid\\ Av. Juan de Herrera 4\\ 28040 Madrid\\ Spain}
\email{daniel.fox@upm.es}
\thanks{D.J.F. has been partially supported by Project PID2024-156578NB-I00 funded by MICIU /AEI /10.13039/501100011033 / FEDER, EU}

\author[V.G.~Tkachev]{Vladimir G. Tkachev}
\address[V.G.~Tkachev]{Department of Mathematics\\ Link\"oping University\\ 58183\\ Sweden}
\email{vladimir.tkatjev@liu.se}
\thanks{V.G.T. has been supported by the Stiftelsen L\"angmanska kulturfonden, Grant  BA26-2924.}


\maketitle

\begin{abstract}
This paper introduces and studies $(2,3)$-palintropic algebras, a class of commutative
algebras defined by the identity $(x^{2})^3 - (x^{3})^2 = 0$. This specific
relation is the simplest generator of the $2$-dimensional space of minimal-degree
evanescent identities in degree $6$, and encompasses several well-studied
structures, including Jordan and medial algebras. The primary motivation for
investigating these algebras lies in their trivial Peirce polynomials, which
removes a priori restrictions on the spectrum of the multiplication operator
associated with an idempotent. This paper reviews and further develops the theory of Peirce operators, Peirce polynomials, and second-order linearizations. It shows that
despite the triviality of the Peirce polynomial, any idempotent $c$ admits
well-behaved, explicit fusion rules for multiplication between its $\lambda$-Peirce
spaces for $\lambda \neq \tfrac{1}{2}$. Furthermore, multiplication
by a nonsingular idempotent is an algebra homomorphism. Finally, there are given concrete examples of $(2,3)$-palintropic algebras, together with 
applications of these algebraic structures to commutative polynomial maps.
\end{abstract}

\tableofcontents

\section{Introduction and motivations}

The usual approach to studying a class of algebras defined by a polynomial identity $P =0$ is to use the Peirce decomposition and corresponding fusion rules associated with an idempotent in an algebra satisfying $P$. These arise from linearization of the defining identity $P$. This approach has been used successfully to study and classify associative algebras, power associative algebras, Jordan algebras, pseudocomposition algebras, Bernstein genetic algebras, Hsiang algebras \cite{Fox-Tkachev}, and many other classes of algebras satisfying identities of low degrees. In \cite{Tkachev-universality}, the case of a general univariate identity (in one formal variable) was considered with the goal of identifying spectral properties of idempotents in algebras satisfying an identity $P$ that depend only on the specific identity $P$.

Linearization yields the \emph{Peirce operator}, $\rho(P;q)$, a $\field$-linear map from the $\field$-algebra of univariate nonassociative commutative polynomials in $x$ (identities) to the polynomial ring $\field[q]$ (see Definition \ref{def:peircepolynomial}). The algebraic set
\begin{align*}
\sigma_{\field}(P) = \{\alpha\in \field: \rho(P; \alpha) = 0\}
\end{align*}
constitutes the \emph{abstract Peirce spectrum} over $\field$ of the identity $P$ (see Definition \ref{def:abstractpeircespectrum}). The Peirce spectrum of any idempotent $c$ contained in a commutative $\field$-algebra satisfying $P$ is contained in the abstract Peirce spectrum of $P$, although the containment can be proper. In this sense the abstract Peirce spectrum of $P$ provides a priori constraints on the Peirce spectra of idempotents in commutative algebras satisfying $P$.

However, as was observed in \cite[p. 230]{Osborn}, it can occur that $\rho(P; q)$ is identically zero. In this case the abstract Peirce spectrum is all of $\field$, $\sigma_{\field}(P) = \field$. In \cite{Tkachev-universality}, an identity was called \emph{degenerated} if its Peirce polynomial $\rho(P; q)$ vanishes identically. The same notion was studied by Varro \cite{varro2020peirce} in a study of baric algebras, where such identities were called \emph{evanescent}, and this terminology is adopted here.

Varro \cite[Proposition $30$]{varro2020peirce} (see Lemma \ref{lem:evanideal}) showed the set $\evan \subset \Nring{\field}{\{x\}}$ of evanescent identities is an ideal and so must be the kernel of a homomorphism to some algebra.
Lemma \ref{lem:starmedial} shows that the polynomial ring $\field[q]$ with the product
$$f\star g = q(f(\tfrac{1}{2})g + g(\tfrac{1}{2})f)
$$ is an (infinite-dimensional) medial algebra in the sense defined in \cite{Krasnov-Tkachev-medial, Tka23a} and recalled below and Theorem \ref{thm:peircehomomorphism} reinterprets the Peirce operator as a homomorphism from $\Nring{\field}{\{x\}}$ to $(\field[q], \star)$.

By their definition, evanescent identities yield no constraints on the values of Peirce eigenvalues of idempotents in algebras satisfying $P$, so such an identity can be said to be \emph{spectrally unconstrained}.
For evanescent identities, the standard approach using the Peirce decomposition encounters difficulties, and it is necessary to revisit the basic methodology. One of the main points of what follows is nonetheless, the Peirce spectrum of a commutative algebra satisfying an evanescent identity $P$ must satisfy multiplicative fusion rules determined by the second-order linearization of the identity.

Let $\evan_{k}$ be the vector space of degree $k$ homogeneous evanescent identities. Varro proved \cite{varro2020peirce} that the smallest degree $k$ for which $\evan_{k}$ is nonzero is $k = 6$ and $\evan_{6}$ has dimension $2$ and is spanned by identities $E_{1}$ and $E_{2}$ defined in \eqref{strange} and \eqref{Sprimeekv} (Proposition \ref{pro:adiabatic}). The identity $E_{1}$ is the main focus of this article in part because it is in some sense the simplest univariate evanescent identity. As is explained next, the $1$-dimensional subspace of $\evan_6$ spanned by $E_{1}$ is distinguished by properties of the second-order linearization of $P$.

Fusion rules for Peirce eigenspaces of an idempotent in an algebra satisfying $P$ are controlled by the second order linearization of $P$. The second-order Peirce polynomial $\frak{D}(P; q, r, s) \in \field[q, r, s]$ is defined from the second order linearization in a manner analogous to how the Peirce operator $\rho(P; q)$ is defined from the first order linearization (see Definition \ref{def:secondorderpeirce}). If a commutative $\field$-algebra $\alg$ satisfies $P$ and $\alpha, \lambda, \mu \in \fie$ are eigenvalues of an idempotent $c \in \alg$ such that the product of $\lambda$ and $\mu$ eigenspaces is contained in a sum of eigenspaces including the $\alpha$-eigenspace, then $\alpha, \lambda, \mu$ must be roots of $\frak{D}(P; q, r, s)$. The \emph{abstract fusion rules} of $P$ are given by the algebraic set
\begin{align*}
\fusion_{\field}(P)
& = \{(\alpha, \lambda, \mu) \in \sigma_{\field}(P)^{3}:  \frak{D}(P; \alpha, \lambda, \mu ) = 0\}.
\end{align*}
(See Definition \ref{def:abstractfusionrules}.) By definition, the fusion rules associated with any idempotent $c$ in a particular algebra satisfying $P$ are constrained by the abstract fusion rules of $P$ in a manner exactly analogous to how the abstract Peirce spectrum of $P$ constrains the Peirce spectrum of $c$.

It can happen that $P$ has the \emph{second-order evanescence property} $\fusion_{\field}(P) = \field^{3}$. Example \ref{ex:secondorderevanescence} gives an example of a homogeneous degree $11$ identity $P$ that is both evanescent and satisfies $\fusion_{\field}(P) = \field^{3}$. However, it turns out that no homogeneous degree $6$ evanescent identity is second-order evanescent in this sense, because $\fusion_{\field}(P)$ is a proper subset of $\field^{3}$ for all such $P$. Moreover, in this case $\frak{D}(P; \alpha, \lambda, \mu )$ is always at most quadratic in $\alpha$ and it is linear in $\alpha$ exactly when $P$ is a multiple of $E_1$ (Proposition \ref{pro:de1e2}). As is explained in more detail later in this introduction, part of the paper is devoted to initiating the study of algebras that satisfy this identity $E_{1}$ and to motivating their study from other perspectives with the goal of better understanding how abstract fusion rules constrain behavior even for evanescent identities.

A general nonassociative commutative monomial in an indeterminate $x$ is denoted by $x^\alpha$, in which the \emph{monomial power} $\alpha$ is a formal symbol that indicates the abstract (non-planar) binary rooted tree corresponding with the monomial $x^{\alpha}$ and coding the manner of parenthesizing the products appearing in $x$. The weight or degree of $\alpha$ means the number of leaves of the corresponding abstract binary rooted tree. The tallest binary rooted trees of a given degree are the  \emph{principal powers} of $x$ are defined inductively by $x^{1} = x$ and $x^{n+1} = x(x^{n})$ for $n \geq 1$.

The subspace of evanescent identities is surprisingly large. In the Peirce polynomial $\rho(x^\alpha; q)$ of the monomial $x^\alpha$ the coefficient of $q^k$ is the number of leaves of height $k$ over the root; see \eqref{peircetreeformula} of Lemma \ref{peirceoperatorlemma}. In particular no monomial is evanescent, but the difference of monomials is evanescent if the corresponding trees have the property that the number of leaves of a given height in either tree is the same for all heights. This observation yields perhaps the simplest examples of evanescent identities,
\begin{equation}\label{strangealphabeta}
(x^{\alpha})^\beta - (x^{\beta})^\alpha =0,
\end{equation}
for monomial powers $\alpha$ and $\beta$ of degree at least $2$ (this is a special case of Proposition \ref{pro:comp}). A commutative algebra $\alg$ is \emph{$(\alpha, \beta)$-palintropic} if it satisfies \eqref{strangealphabeta}. Given two integers $2\le p<q$, an algebra $\alg$ is $(p,q)$-\emph{palintropic} if any $x\in\alg$ satisfies
$(x^{p})^q=(x^{q})^p$.
The simplest homogeneous evanescent identity is the previously mentioned degree $6$ identity
\begin{equation}\label{strange}
E_1:=(x^{2})^3 - (x^{3})^2=0,
\end{equation}
that is the particular case $(p,q)=(2,3)$.
The medial algebras studied in \cite{Krasnov-Tkachev-medial, Tka23a} satisfy $E_1$, and this is another motivation for studying $(2, 3)$-palintropic algebras.

\begin{remark}
The terminology originates in \cite{EtheringtonSomeNAAmulti,Etherington1945transposed} where an algebra is called \emph{palintropic} if $(x^{p})^q=(x^{q})^p$ for \emph{all} principal powers.
Etherington remarks that  unital algebras of degree two and Jordan algebras are palintropic.
In \cite{EtheringtonSomeNAAmulti}, Etherington also gives another example of a palintropic algebra $\M_n$ defined by the multiplication rules $a_ia_j=0$ for $i \neq j$ and $a_ia_i=a_{i+1}$ (indices taken cyclically) in a basis $\{a_1,\ldots,a_n\}$.
Here it serves to show that there exist nontrivial $(p, q)$-palintropic algebras.
The specific identity \eqref{strange} appears in several contexts, for example in \cite{EtheringtonSomeNAAmulti} and in \cite[sec.~3]{mallol2017varro} in the context of backcrossing algebras \cite{mallol1994varro}. \qed
\end{remark}

Next the methodology of the Peirce decomposition is recalled.
Here $\alg=(\alg,\mult)$ denotes a commutative (not necessarily associative) algebra over a field $\field$ of characteristic not $2$ with multiplication $\mult$. The assumption $\operatorname{char}\field \neq 2$ is supposed implicitly throughout the paper because $\tfrac{1}{2}$ appears constantly. From Section \ref{sec:degreesix} onwards, it is also assumed that $|\field| > 6$.
Although $\alg$ is allowed to be infinite-dimensional, when it is finite-dimensional its dimension will be denoted $\dim \alg=n$.
Throughout the paper, when the algebra product $\mult$ is clear from context, the notation indicating it is omitted and the product of elements $x, y\in \alg$ is indicated by juxtaposition, $xy$, instead of $x\mult y$. Notation such as $\mult$ indicating explicitly the product is used when it is helpful to distinguish various different products under consideration.
The operator $L(a):\alg \to \alg$ of multiplication by $a$ is defined by $L(a)x = ax$.
An element $c\in \alg$ is an idempotent if $c^2=c$. The set of nonzero idempotents of $\alg$ is denoted $\Idm(\alg)$.

If $\alg$ is finite-dimensional, for $a \in \alg$ the set of roots of $\det(L(a)-\lambda\Id)=0$ in an algebraic closure $\bar\field$ of $\field$ is denoted $\sigma(a)$ and the \emph{Peirce spectrum} of the algebra $\alg$ is defined to be
$$
\sigma(\alg)= \bigcup_{c\in \Idm(\alg)} \sigma(c).
$$
As the existence of idempotents depends on the arithmetic properties of the underlying field, the spectrum $\sigma(\alg)$ is sensitive to the underlying field. For $c \in \Idm(\alg)$, the eigenspace $\alg_c(\lambda) = \ker(L(c) - \lambda \Id)$ is called the \emph{Peirce subspace} associated with $\lambda \in \sigma(c)$.

If a $\field$-algebra is given by an identity $P$, then the Peirce spectrum of $\alg$ is contained in the abstract Peirce spectrum of $P$,
$$\sigma(\alg)\cap \field \subset \sigma_{\field}(P),$$
so that the identity imposes a priori restrictions on the Peirce spectrum $\sigma(\alg)$ provided $P$ is not evanescent.

Because the degree of the Peirce polynomial $\rho(P; q)$ of a homogeneous identity $P$ is less than or equal to $\deg P - 1$ (Corollary \ref{cor:degpeirce}), if the Peirce polynomial $\rho(P;q)$ is nonzero, the cardinality of the abstract Peirce spectrum of $P$ is at most $\deg P - 1$. In particular, when a commutative algebra satisfies an identity $P$ that is not evanescent, its Peirce spectrum is finite, with cardinality bounded by $\deg P - 1$, although in general it can be still smaller (the concept of \emph{Peirce complexity} introduced in \cite{FoxTka26} makes precise the possible gap). This is in sharp contrast to the situation for evanescent identities, for which the cardinality of abstract Peirce spectrum is unbounded.

Practically, the abstract Peirce spectrum is extracted from the defining algebraic identity by linearization followed by specialization to an algebra idempotent. Precisely, replacing $x$ in the linearization of $P(x)$ by an idempotent $c$ and simplifying the resulting expression yields the Peirce polynomial evaluated at $q = L(c)$; see Proposition~\ref{pro:Peircedef} for the precise statement.
Example \ref{ex:1} illustrates this process in one of the simplest cases, that of power associative algebras.

\begin{example}
\label{ex:1}
An algebra $\alg$ is power-associative if the subalgebra generated by any element $x$ is associative.
By \cite[Theorem~1]{Albert-powerassociative}, if $\operatorname{char}\field$ is prime to 30 then a commutative algebra $\alg$ is power-associative if and only if it satisfies the degree $4$ identity
\begin{equation}\label{pow}
A(x):=x^4-x^{2}x^2 = 0.
\end{equation}
In this case, linearization of \eqref{pow} yields
\begin{equation}\label{Lpow}
2L(x)^3+L(x)L(x^2)+L(x^3) -4L(x^2)L(x)=0,
\end{equation}
and conversely, \eqref{Lpow} implies \eqref{pow}.
For $x$ equal to an idempotent $c\in \Idm(\alg)$, \eqref{Lpow} becomes
\begin{equation}\label{Lpowc}
(2L(c)-1)(L(c)-1)L(c) = 2L(c)^3-3L(c)^2+L(c) = 0.
\end{equation}
The polynomial in $L(c)$ appearing in \eqref{Lpowc},
$$\rho(A; q) = 2q^{3} - 3q^{2} + q,$$
is the Peirce polynomial of the identity $A$, and the abstract Peirce spectrum of $A$ comprises its roots, $\sigma_{\field}(A)=\{0,\tfrac12,1\}$.
An associative algebra provides a simple example of a power associative algebra for which the Peirce spectrum $\sigma(c)$ of any nonzero idempotent is contained properly in the abstract Peirce spectrum, for $\tfrac{1}{2}$ is not in the spectrum of an idempotent of an associative algebra.

If $L(c)$ is moreover diagonalizable, then there holds the corresponding Peirce decomposition $\alg=\alg_c(1)\oplus\alg_c(0)\oplus\alg_c(\tfrac12)$.
Utilizing the second linearization of \eqref{Lpow} yields the \emph{fusion rules} describing the containments of products of Peirce eigenspaces,
\begin{align}\label{fusionjordan}
&\alg_c(\lambda) \alg_c(\lambda)\subset \alg_c(\lambda),&
&\alg_c(\lambda) \alg_c(\tfrac12)\subset \alg_c(\tfrac12),&
&\alg_c(\tfrac12)\alg_c(\tfrac12)\subset \alg_c(0)\oplus \alg_c(1),
\end{align}
for $\lambda \in \{0, 1\}$ (see \cite[Theorem~2]{Albert-powerassociative}.)
\qed\end{example}

The study of algebras satisfying irreducible identities of degree greater than $4$ is somewhat limited, in part because of the complexity of the Peirce decomposition methodology, although baric algebras and train algebras are some well known exceptions. One example is the algebras satisfying the degree $5$ identity
\begin{equation}\label{osborn5}
2x^5-3(x^2)^2x+x^3x^2=0,
\end{equation}
introduced by Osborn in
\cite{osborn1965identities}.
The abstract Peirce spectrum of \eqref{osborn5} is $\{0,1,1/2\}$; see Example~\ref{ex:osborn} below. This makes the structure of the Osborn algebras similar to that of Jordan and power-associative algebras, which they generalize, although the fusion rules for Osborn algebras are more complicated.
In this case the eigenvalue $\tfrac{1}{2}$ has multiplicity $2$ as a root of the Peirce polynomial. Multiplicity is an interesting aspect not captured by the definition of the abstract Peirce spectrum as a zero locus.

In practice, a focus on low degree identities hides the phenomenon of evanescence, as the lowest degree of a homogeneous univariate evanescent identity is $6$.

A different motivation for studying $(2,3)$-palintropic algebras is an analogy relating associative and medial algebras in which $(2,3)$-palintropic algebras are analogues of power associative algebras.

In a commutative algebra the associator $[x, y, z] = (xy)z - x(yz)$ measures the failure of equality of the products corresponding with different labelings of the abstract binary rooted tree corresponding with the third principal power. Similarly, the mediator of Definition \ref{defmediator} measures the failure of equality of the products corresponding with different labelings of the abstract binary rooted tree corresponding with the monomial $(x^2)^2$.

\begin{definition}\label{defmediator}
Let $\alg$ be a commutative $\field$-algebra (possibly infinite-dimensional). The \emph{mediator} is the $\field$-multilinear map from $\otimes^{4}\alg$ to $\alg$ defined for $x,y,z,w \in \alg$ by
\begin{align}
[x, y, z, w]:=(xy)(zw) - (xz)(yw).
\end{align}
\end{definition}

The vanishing of the associator defines the class of associative algebras. The vanishing of the mediator yields the following class of algebras. A commutative algebra $\alg$ is \emph{medial} \cite{Krasnov-Tkachev-medial, Tka23a} if it satisfies the identity
\begin{equation}\label{medialdef}
[x,y,z,w]=0.
\end{equation}

Consider evaluating the associator or mediator on monomials in $x$ in a commutative algebra. The lowest degree case for which the associator of various monomials in $x$ is not automatically zero is $[x, x, x^{2}]$, which equals $(x^2)^2 - x^4$ and whose vanishing is equivalent to power associativity if the characteristic is prime to $30$. Similarly, the lowest degree case for which the mediator of monomials in $x$ is not automatically zero is $[x, x, x^{2}, x^{2}]$, which equals $(x^2)^3 - (x^3)^2$ and whose vanishing characterizes the $(2, 3)$-palintropic algebras. Thus treating the associator and mediator in formally the same way yields a parallelism between associativity and mediality and between power associativity and $(2, 3)$-palintropicity.

A further connection is seen by considering unital commutative algebras. Theorem~5.5 in \cite{Tka23a} states that the category of special medial algebras is isomorphic to the category of unital calibrated associative algebras.
A bit more prosaically, evaluating the mediator on a unit $e$ shows that $[x, y, z, e] = [y, x, z]$. In particular, unital medial algebras are associative. Proposition \ref{pro:unital} shows that, analogously, unital $(2,3)$-palintropic algebras are power associative.

The preceding discussion shows several senses in which $(2, 3)$-palintropic algebras are a natural generalization of power associative algebras. Consideration of examples further confirms this perspective.

A commutative algebra $\alg$ is Jordan if it satisfies
\begin{equation}\label{jordandef}
J(x,y):=x^2(xy)-x(x^2y)=0.
\end{equation}
Define two superclasses of medial algebras, the first of which is a subclass of the other. An algebra is \emph{J-medial} (Jordan medial) if it satisfies
\begin{equation}\label{medialJ}
K(x,y,z):=[x,y,z,x^2]=(xy)(x^2z)-(x^2y)(xz)=0.
\end{equation}
The J-medial algebras are a subclass of those satisfying the identity
\begin{equation}\label{pseudoJ}
P(x,y):=K(x,y, x^2)=[x,y, x^2,x^2]=(x^2)^2(xy) -x^3(x^2y) = 0.
\end{equation}
Since
\begin{equation}\label{Siden}
[x,x, x^2,x^2]=K(x,x, x^2)=P(x,x)=E_1(x),
\end{equation}
there result the following inclusions:
\begin{equation}\label{inclus}
\boxed{\text{Medial}}\subset\boxed{\text{$J$-Medial}}\subset\boxed{\text{satisfying \eqref{pseudoJ}}}\subset\boxed{\text{$(2,3)$-palintropic}}
\end{equation}
There is another chain of inclusions (the second inclusion is valid for the ground field of characteristic not 2):
\begin{equation}\label{inclus2}
\boxed{\text{Associative}}
\subset\boxed{\text{Jordan}}
\subset
\boxed{\text{Power-associative}}
\subset\boxed{\text{$(2,3)$-palintropic}}
\end{equation}
Although there are algebras that belong to classes in both chains, the individual classes in \eqref{inclus} and \eqref{inclus2} are all pairwise distinct.
See, however, Proposition~\ref{pro:pseudounital} below which shows that any \textit{unital} algebra satisfying \eqref{pseudoJ} is Jordan. Note also that any \emph{unital} $J$-medial algebra is (unital) Jordan: substitution of $z=e$ into \eqref{medialJ} yields
$$
K(x,y,e)=(xy)(x^2e)-(x^2y)(xe)=(xy)x^2-(x^2y)x=J(x,y).
$$
This is parallel to the relations between unital medial and unital $(2,3)$-palintropic algebras and associative and power associative algebras.

Proposition \ref{pro:classes} shows that the $(2,3)$-palintropic algebras constitute a large class that includes a wide variety of algebras, among them Jordan, power-associative and medial algebras whose classifications are already rather nontrivial.

\begin{proposition}\label{pro:classes}
The following algebras are $(2, 3)$-palintropic:
\begin{enumerate}
[label=$\mathrm{(\roman*)}$,itemsep=0.5ex,leftmargin=1cm]
\item\label{i:tri} any commutative nilalgebra $\alg$ of nilindex $3$ ($x^{3} =0$ for all $x \in \alg$);
\item\label{i:pow} any commutative power-associative algebra;
\item\label{i:jord} any Jordan algebra;
\item\label{i:med} any medial algebra;
\item\label{i:Jmed} any J-medial algebra;
\item\label{i:pseud} any commutative algebra satisfying \eqref{pseudoJ}.
\end{enumerate}
\end{proposition}

\begin{proof}
Claims \ref{i:tri} and \ref{i:pow} are immediate.
Because the characteristic is not $2$, a Jordan algebra is power associative, so case \ref{i:jord} follows from case \ref{i:pow}. By \eqref{Siden}, \ref{i:med}, \ref{i:Jmed}, and \ref{i:pseud} are $(2, 3)$-palintropic.
\end{proof}

Example \ref{ex:3} exhibits a $(2, 3)$-palintropic algebra that is not medial.
It is not clear how large is the class of $(2, 3)$-palintropic algebras not included in Proposition~\ref{pro:classes} and it is unclear what form a classification theory for general $(2,3)$-palintropic algebras might take.
It is apparent that a direct sum of $(2,3)$-palintropic algebras is again a $(2,3)$-palintropic algebra, and so it is reasonable to expect that classification can be reduced to consideration of the simple and nilpotent cases.
Moreover, it seems plausible that the classification of simple $(2,3)$-palintropic algebras is tractable, although this is not addressed here.

The Peirce spectrum of a $(2, 3)$-palintropic algebra can be quite varied.
A nilalgebra as in \ref{i:tri} of Proposition~\ref{pro:classes} contains no nonzero idempotents (over any extension of the ground field).
If a power-associative algebra has an idempotent, its Peirce spectrum is $\{0,\frac12,1\}$.
Medial algebras in general can have rather varied spectra, see \cite{Tka23a} and also Proposition~\ref{pro:severalvar} below.

For commutative algebras satisfying a univariate identity, a special role is played by the eigenvalue $\tfrac{1}{2}$ \cite{Tkachev-universality}.
\begin{definition}
Let $\alg$ be a commutative $\field$-algebra. An idempotent $c \in \Idm(\alg)$ is \emph{singular} if the endomorphism
\begin{align}\label{Phidefined}
\Phi(c): = \Id - 2L(c)
\end{align}
is not invertible and $c$ is \emph{nonsingular} otherwise.
\end{definition}
If $\alg$ is finite dimensional, $c$ is singular if and only if $\tfrac{1}{2} \in \sigma(c)$.
For any idempotent $c$ in a medial algebra $\alg$, $L(c)$ is an algebra homomorphism; this follows from taking  $z = w = c$ in \eqref{medialdef}. Theorem \ref{the:iii} shows that for any nonsingular idempotent $c$ in a $(2, 3)$-palintropic algebra, $L(c)$ is an algebra homomorphism, so that this property of medial algebras extends partially to all $(2, 3)$-palintropic algebras.

In \cite{Tka23a} it was shown that any medial algebra that is \emph{special}, meaning it contains an \emph{invertible} idempotent $c$ (meaning $L(c)$ is invertible) is isotopic, via $L(c)^{-1}$, to an associative algebra. It would be interesting to extend this result to $(2, 3)$-palintropic algebras. Proposition \ref{pro:unital} shows that a unital $(2, 3)$-palintropic algebra is power associative and Corollary \ref{cor:subpa} shows the Peirce subspace $\alg_{c}(1)$ of a nonzero idempotent $c$ is a power associative subalgebra. This makes it reasonable to ask: when is a $(2, 3)$-palintropic algebra containing a nonsingular invertible idempotent isotopic to a power associative algebra?

The final Section \ref{sec:integrable} explains that certain evanescent identities such as the $(p, q)$-palintropic identity give rise to polynomial mappings of $\field^{n}$ that are integrable in the sense defined by Veselov \cite[Chap.~2]{veselov1991integrable}.
The $(p, q)$-palintropic identity can be thought of as the commutativity relation $\phi_p\circ \phi_q = \phi_q\circ \phi_p$
where $\phi_k(x)=x^k$. This gives a general construction of integrable polynomial mappings whose relation with previously known constructions of such mappings remains to be clarified.

\section{Peirce operator and Peirce polynomial}\label{sec:Peirce}
This section defines the Peirce operator and the Peirce polynomial of a nonassociative polynomial.

It is convenient to begin summarizing notation and material relating polynomial identities, nonassociative monomials, and abstract binary rooted trees. This follows loosely the second author's \cite{Tkachev-universality}, as well as \cite{Osborn} and \cite{Holtkamp, VanDerLinden, Zhevlakov}, adapted to the circumstance that here there are considered only commutative algebras.

Given a finite set $\mathsf{T}$, $\nass(\mathsf{T})$ denotes the free nonassociative commutative magma generated by $\mathsf{T}$.
Its elements are finite nonassociative commutative parenthesized words in the elements of $\mathsf{T}$ and its product is the commutative concatenation of such parenthesized words. It can be identified with the set of abstract (non-planar) binary rooted trees whose leaves are labeled by elements of $\mathsf{T}$ equipped with the product given by grafting of rooted trees \cite{Holtkamp,Petrogradsky-trees,  Tkachev-universality, varro2020peirce} (more details are given below and in the proof of Lemma \ref{peirceoperatorlemma}). Each $x \in \mathsf{T}$ determines a word $(x)$ customarily denoted simply by $x$, with the parentheses omitted. More generally, the outermost parentheses of a word are omitted for readability. If $\mathsf{T} = \{x, y, z\}$, some typical words are $x$, $y$, $xy$, $x(yz)$, and $(xy)(xz)$. The commutativity means that two words are identified if they correspond with the same abstract binary rooted tree, so that one can be obtained from the other by permuting subwords within each reciprocal pair of parentheses. For example, the words $x(yz)$, $(yz)x$, and $(zy)x$ are all identified.

Given a field $\field$, $\Nring{\field}{\mathsf{T}}$ denotes the free commutative $\field$-algebra generated by $\nass(\mathsf{T})$.
Its elements are finite linear combinations of elements of $\nass(\mathsf{T})$ and its product is given by extending the commutative magma product on $\nass(\mathsf{T})$ bilinearly over $\field$.

An element of $\Nring{\field}{\mathsf{T}}$ is called a \emph{nonassociative polynomial} over $\field$. A nonassociative polynomial is called a \emph{nonassociative monomial} if it is a scalar multiple of an element of $\nass(\mathsf{T})$.

The \emph{degree} $\weight_{x}(t)$ of $x \in \mathsf{T}$ in the word $t \in \nass(\mathsf{T})$ is the number of times $x$ appears in $t$. The degree $\deg_{x}(P)$ of $P \in \Nring{\field}{\mathsf{T}}$ is the maximum of the degree of the words in an expression of $P$ as a linear combination of elements of $\nass(\mathsf{T})$. The \emph{total degree} $\deg(P)$ of $P \in \Nring{\field}{\mathsf{T}}$ is defined by $\deg(P)=\sum_{x \in\mathsf{T}}\deg_{x}(P)$. For example, $\weight_{t_1}((t_1t_1)(t_2t_3))=2$ and $\deg ((t_1t_1)(t_2t_3))=4$. The degree in $x$ and the total degree are additive with respect to the product in $\Nring{\field}{\mathsf{T}}$, meaning $\deg_x(PQ) = \deg_x(P) + \deg_x(Q)$ and similarly for $\deg$.

In the simplest case, $\mathsf{T}=\{x\}$, there is a bijection between $\nass(\mathsf{\{x\}})$ and abstract binary rooted trees with leaves labeled by $x$. An abstract binary rooted tree $T$ means a connected graph without cycles having a distinguished vertex, called the root, and such that the valence of the root is $0$ or $2$ and the valence of any vertex different from the root is $1$ or $3$. Its leaves $L(T)$ are the vertices of valence $1$ (and the root when it has valence $0$). The \emph{graft} $T_{1}\vee T_{2}$ of binary rooted trees $T_{1}$ and $T_{2}$ is formed by connecting a new vertex, the root of the graft, to the root vertex of each of $T_{1}$ and $T_{2}$. This makes the set of isomorphism classes of abstract binary rooted trees a commutative magma isomorphic with $\nass(\mathsf{\{x\}})$. A generic nonassociative monomial is denoted by $x^\alpha$, in which $\alpha$ is a formal symbol that indicates the abstract binary rooted tree corresponding with the monomial $x^{\alpha}$. The isomorphism associates with $T$ the monomial $x^{\alpha}$ given by parenthesizing the labels of leaves adjacent to a common vertex. The degree of the monomial $x^\alpha$ equals the number of leaves of the corresponding tree.

For the subspace $\Nring{\field}{\{x\}}_n \subset \Nring{\field}{\{x\}}$ spanned by monomials of degree $n$, the dimension
\begin{equation}\label{wedder}
\dim \Nring{\field}{\{x\}}_n=W_n,
\end{equation}
is the Wedderburn-Etherington number $W_n$ \cite{etherington1937non, wedderburn1922functional} that counts the number of abstract binary rooted trees with $n$ leaves and whose smallest values are:

\begin{center}
\begin{tabular}{c|r|r|r|r|r|r|r|r|r|c}
  $n$ & $1$ & $2$ & $3$ & $4$ & $5$ & $6$ & $7$ & $8$ &9& $10$ \\\hline
  $W_n$ & $1$ & $1$ & $1$ & $2$ & $3$ &$6$ & $11$ & $23$ &46& $98$ \\
\end{tabular}
\end{center}

If $v$ is a vertex of the tree $T$, the height $\height(v)$ is the length of the shortest path from $v$ to the root vertex. The height of the tree $T$ is the maximum of the heights of its leaves.
The \emph{principal powers} of $x$ are defined inductively by $x^{1} = x$ and $x^{n+1} = x(x^{n})$ for $n \geq 1$. The height of the tree $T$ is at most one less than the degree of the corresponding monomial and the principal power $x^{n+1}$ corresponds with the unique abstract binary rooted tree of degree $n+1$ having height $n$.

The unordered multiset $\height(T) := \{\height(v): v \in L(T)\}$ comprises the heights of the leaves of $T$.
That $\height(T)$ satisfies
\begin{align}\label{binary1}
\sum_{v \in L(T)}2^{-\height(v)} = 1
\end{align}
can be proved by induction on the number of vertices of $T$ (the inductive step is the $q = 1/2$ specialization of \eqref{peirceheight} below). The essential content of Lemma \ref{treelemma} is that any set of putative heights satisfying the condition \eqref{binary1} can be realized as the heights of some abstract rooted binary tree.

\begin{lemma}\label{treelemma}
The following sets are in bijection:
\begin{enumerate}
[label=$\mathrm{(\mathbf{NC}\arabic*)}$,
itemsep=0.5ex,
leftmargin=3.5em
]
\item\label{NA:1}
nonassociative monomials $x^\alpha\in \nass(\{x\})$ of degree $m$;
\item\label{NA:4}
abstract binary rooted trees with $m$ leaves,
\item\label{NA:5}
Unordered multisets $\{k_{1}, \dots, k_{m}\}$ of nonnegative integers such that $\sum_{i = 1}^{m}2^{-k_{i}} = 1$.
\end{enumerate}
\end{lemma}

\begin{proof}
The bijection between \ref{NA:1} and \ref{NA:4} is well known and has been sketched above. Essentially an abstract binary rooted tree with leaves labeled by $x$ is notation indicating parenthesizations in the corresponding monomial $x^{\alpha}$.
That an unordered multiset $\{k_{1}, \dots, k_{m}\}$ of nonnegative integers equals the set $\height(T)$ of heights of an abstract binary rooted tree if and only if $\sum_{i = 1}^{m}2^{-k_{i}} = 1$ follows from Exercise~$3$ in Section~2.3.4.5 of \cite{knuth1997art} (and its solution on p. $594$ of \cite{knuth1997art}).
\end{proof}

Lemma \ref{peirceoperatorlemma} defines the Peirce operator.
\begin{lemma}\label{peirceoperatorlemma}
There is a unique $\field$-linear map $\rho: \Nring{\field}{\{x\}} \to \field[q]$ that associates with $P \in  \Nring{\field}{\{x\}}$ an (associative) polynomial in the indeterminate $q$, $\rho(P; q) \in \field[q]$, and satisfying the initial condition
and recurrence relation
\begin{align}
\label{monomialpeircerecursion}
\begin{aligned}
\mathrm{(i)}\,\,&\rho(x;q)=1,\\
\mathrm{(ii)}\,\,&\rho(x^{\alpha}x^{\beta};q)=q(\rho(x^{\alpha};q)+\rho(x^{\beta};q)),
\end{aligned}
\end{align}
for all $x^{\alpha}, x^{\beta} \in \Nring{\field}{\{x\}}$. Moreover, if $T$ is the tree corresponding with $x^{\alpha}\in  \Nring{\field}{\{x\}}$,
then $\rho(x^{\alpha}; q)$ is given by
\begin{align}\label{peircetreeformula}
\rho(x^{\alpha}; q) =\sum_{v \in L(T)}q^{\height(v)},
\end{align}
where the sum runs over the leaves of $T$, $L(T)$, and $\height(v)$ is the height of the leaf $v$.
\end{lemma}

\begin{proof}
The essential points of the proof are indicated. Consult \cite[Section $4$]{Tkachev-universality} and \cite[Section $3$]{varro2020peirce} for details.
The map $\rho$ is defined on $\Nring{\field}{\{x\}}$ by extending \eqref{peircetreeformula} $\field$-linearly. It is straightforward to check that this satisfies \eqref{monomialpeircerecursion}. In particular, the recursion relation \eqref{monomialpeircerecursion} follows from the observation that if $T = T_{1}\vee T_{2}$, then
\begin{align}\label{peirceheight}
\begin{aligned}
\sum_{v \in L(T)}q^{\height(v)} &= \sum_{v \in L(T_{1})}q^{\height(v)} + \sum_{v \in L(T_{2})}q^{\height(v)} = q\Big(\sum_{v \in L(T_{1})}q^{\height_{1}(v)} + \sum_{v \in L(T_{2})}q^{\height_{2}(v)}\Big),
\end{aligned}
\end{align}
where the height $\height_{i}(v)$ of $v\in T_{i}$ in $T_{i}$ satisfies $\height(v) = \height_{i}(v) + 1$.
\end{proof}

\begin{corollary}\label{cor:degpeirce}
$\deg_q \rho(P;q)\le \deg_{x}(P)-1$.
\end{corollary}
\begin{proof}
This follows from \eqref{peircetreeformula} because the height of a leaf in the tree corresponding with a monomial is at most the degree of the monomial minus $1$.
\end{proof}

\begin{definition}\label{def:peircepolynomial}
Given $P\in \Nring{\field}{\{x\}}$, $\rho(P;q)$ is called the \emph{Peirce polynomial} of $P$. The operator $\rho$ of Lemma \ref{peirceoperatorlemma} is the \emph{Peirce operator}.
\end{definition}

\begin{remark}
The Peirce operator sends a nonassociative polynomial to an ordinary associative polynomial. It can be thought of as an \emph{associativization} map determined by counting parentheses. \qed
\end{remark}

\begin{example}\label{principalpowersexample}
Induction on $n$ using \eqref{monomialpeircerecursion} shows that, for $n \geq 2$, the Peirce polynomial of the $n$th principal power of $x$ is
\begin{align}\label{peircexn}
\rho(x^n; q) = 2q^{n-1} + \sum_{j = 1}^{n-2}q^{j}.
\end{align}
\qed\end{example}

For $P \in \Nring{\field}{\{x\}}$, $P(1)$ denotes the image of $P$ under the evaluation homomorphism $\Nring{\field}{\{x\}} \to \fie$ given by evaluating $x$ at $1 \in \fie$.
\begin{lemma}\label{peircerecursionlemma}
For all $P, Q \in \Nring{\field}{\{x\}}$ there holds
\begin{align}
\label{peircerecursion}
\rho(PQ;q)&=q(Q(1)\rho(P;q)+P(1)\rho(Q;q)).
\end{align}
\end{lemma}
\begin{proof}
This follows from \eqref{monomialpeircerecursion} by $\fie$-linearity of the Peirce operator.
\end{proof}

\begin{example}
Inductions on the degree of $\alpha$ using \eqref{monomialpeircerecursion} show that, for a monomial $x^\alpha\in \Nring{\field}{\{x\}}$,
\begin{align}
\label{PEIRCE3}&\rho(x^\alpha;1)= \deg x^\alpha,\\
\label{PEIRCE4}&\rho(x^\alpha;\half12)=1,\\
\label{PEIRCE5}&\rho(x^\alpha;0)=
\begin{cases}
   1&\text{if}\,\,\deg x^\alpha=1\\
0&\text{if}\,\,\deg x^\alpha\ge 2,
  \end{cases}
\end{align}
in which the first equality makes sense in positive characteristic if $\deg x^{\alpha}$ is viewed as an element of $\field$ and the second equality requires $ \operatorname{char}\field\neq 2$. For principal powers $x^{n}$, the identities \eqref{PEIRCE3}-\eqref{PEIRCE5} result from evaluating \eqref{peircexn} at $q \in \{0, 1/2, 1\}$. The identities \eqref{PEIRCE3}-\eqref{PEIRCE5} provide useful checks for by-hand computations.
\qed\end{example}

\begin{corollary}\label{cor:evaluationP1}
Any $P \in \Nring{\field}{\{x\}}$ satisfies $P(1) = \rho(P; \tfrac{1}{2})$.
\end{corollary}
\begin{proof}
Write $P$ as a linear combination of monomials and use \eqref{PEIRCE4}.
\end{proof}

Corollary \ref{monomialimagecorollary} characterizes the polynomials $Q$ that occur as the Peirce polynomial of a nonassociative power  $x^\alpha$.

\begin{corollary}\label{monomialimagecorollary}
Given a polynomial $Q\in \mathbb{Z}[t]$ with nonnegative integer coefficients,
there exists a nonassociative monomial $x^\alpha$ such that $Q(q)=\rho(x^\alpha;q)$ if and only if $Q(\tfrac12)=1$.
\end{corollary}
\begin{proof}
This follows from \eqref{binary1}, \eqref{PEIRCE4}, and Lemma \ref{treelemma}.
\end{proof}

\begin{example}
The preceding is illustrated with a concrete example. 
The polynomial $Q(t)=6t^6 +3t^5 + 7t^4 +3t^3$ satisfies $Q(\tfrac{1}{2}) = 1$, so, by Corollary \ref{monomialimagecorollary}, $Q(q)$ is the Peirce polynomial of a nonassociative monomial, which can be seen as follows to be the monomial
$$
(x^2x^3)((x^2)^2  (x^3(x^3(x^2)^2)))
$$
associated with the tree of Figure \ref{fig:forest}.
Consider the multiset
$$
A=\{(\tfrac1{2^6})^{\#6},(\tfrac1{2^5})^{\#3}, (\tfrac1{2^4})^{\#7}, (\tfrac1{2^3})^{\#3}\},
$$
where $a^\#$ indicates the multiplicity with which the number $a$ appears. This yields a decomposition of $1$ as a sum of binary powers:
\begin{equation}\label{Adefin}
|A|=\frac{6}{2^6} +\frac{3}{2^5} + \frac{7}{2^4} +\frac{3}{2^3} =1.
\end{equation}
A particular grouping of the summands
\begin{align*}
\underbrace{\underbrace{\textcolor{red}{\frac{1}{2^3}+\frac1{2^3}}}_{=\frac14}
+
\underbrace{\textcolor{blue}{\frac1{2^3}}+\textcolor{blue}{\frac1{2^4}+\frac1{2^4}}}_{=\frac14}}_{=\frac12}
&+
\underbrace{\underbrace{\textcolor{red}{\frac1{2^4}+\frac1{2^4}}+\textcolor{green!60!black}{\frac1{2^4}+\frac1{2^4}}}_{=\frac14}
+
\underbrace{\textcolor{blue}{\frac1{2^4}+\frac1{2^5}+\frac1{2^5}+\frac1{2^5}+\frac1{2^6}+\ldots+\frac1{2^6}}}_{=\frac14}}_{=\frac12}
\end{align*}
corresponds with the binary tree with height sequence $$[3,3,3,4,4,4,4,4,4,4,5,5,5,6,6,6,6,6,6],
$$
as in Figure \ref{fig:forest}. By \eqref{peircetreeformula} this yields
\begin{equation}\label{Knuth2}
\rho(x^\alpha;q)=6q^6 +3q^5 + 7q^4 +3q^3.
\end{equation}
\begin{figure}[ht]
\begin{forest}
  for tree={
  l-=10pt,
    s sep=0.55cm,
    where n children=0{
      circle,
      inner sep=1.5pt,
      thick,
      baseline,
      tier=leaves
    }{
      coordinate
    },
    edge={thick}
  },
  /tikz/red edge/.style={draw=red!70},
  /tikz/blue edge/.style={draw=blue!80!black},
  /tikz/black edge/.style={draw=green!60!black},
  /tikz/red leaf/.style={draw=red!70, label={below:\small\sffamily 3}},
  /tikz/red leaf four/.style={draw=red!70, label={below:\small\sffamily 4}},
  /tikz/black leaf/.style={draw=green!60!black, label={below:\small\sffamily 4}},
  /tikz/blue leaf three/.style={draw=blue!80!black, label={below:\small\sffamily 3}},
  /tikz/blue leaf four/.style={draw=blue!80!black, label={below:\small\sffamily 4}},
  /tikz/blue leaf five/.style={draw=blue!80!black, label={below:\small\sffamily 5}},
  /tikz/blue leaf six/.style={draw=blue!80!black, label={below:\small\sffamily 6}},
  fill=red!70, circle, inner sep=1.5pt,
  [
    [, edge={red edge}
      [, edge={red edge}
        [, edge={red edge}, red leaf]
        [, edge={red edge}, red leaf]
      ]
      [, edge={blue edge}
        [, edge={blue edge}, blue leaf three]
        [, edge={blue edge}
          [, edge={blue edge}, blue leaf four]
          [, edge={blue edge}, blue leaf four]
        ]
      ]
    ]
    [, edge={blue edge}
      [, edge={blue edge}
        [, edge={red edge}
          [, edge={red edge}, red leaf four]
          [, edge={red edge}, red leaf four]
        ]
        [, edge={black edge}
          [, edge={black edge}, black leaf]
          [, edge={black edge}, black leaf]
        ]
      ]
      [, edge={blue edge}
        [, edge={blue edge}
          [, edge={blue edge}, blue leaf four]
          [, edge={blue edge}
            [, edge={blue edge}, blue leaf five]
            [, edge={blue edge}, blue leaf five]
          ]
        ]
        [, edge={blue edge}
          [, edge={blue edge}
            [, edge={blue edge}, blue leaf five]
            [, edge={blue edge}
              [, edge={blue edge}, blue leaf six]
              [, edge={blue edge}, blue leaf six]
            ]
          ]
          [, edge={blue edge}
            [, edge={blue edge}
              [, edge={blue edge}, blue leaf six]
              [, edge={blue edge}, blue leaf six]
            ]
            [, edge={blue edge}
              [, edge={blue edge}, blue leaf six]
              [, edge={blue edge}, blue leaf six]
            ]
          ]
        ]
      ]
    ]
  ]
\end{forest}
\caption{The tree of $(x^2x^3)((x^2)^2(x^3(x^3(x^2)^2)))$}
\label{fig:forest}
\end{figure}\qed
\end{example}

\begin{lemma}\label{lem:starmedial}
On $\field[q]$ define a $\field$-linear map $\tau:\field[q] \to \field$, a commutative multiplication $\star$, and a symmetric bilinear form $\omega$ by
\begin{align*}
\tau(f) &:= f(\tfrac{1}{2}),\\
f\star g &:= q(\tau(f)g + \tau(g)f) = q(f(\tfrac{1}{2})g + g(\tfrac{1}{2})f),\\
\omega(f, g) &:= \tau(f)\tau(g) = f(\tfrac{1}{2})g(\tfrac{1}{2}) = (f\star g)(\tfrac{1}{2}) = \tau(f\star g).
\end{align*}
The algebra $(\field[q], \star)$ is medial with associator satisfying
\begin{align}\label{starassociator}
[f, g, h] &= (q^{2} - q)(\omega(g, h)f - \omega(f, g)h) = (q^{2} - q)\tau(g)\begin{vmatrix} f & h\\ \tau(f) & \tau(h)\end{vmatrix},
\end{align}
for $f, g, h \in \field[q]$.
\end{lemma}

\begin{proof}
The only claim that is not immediate is that $(\field[q], \star)$ is medial. By the definitions of $\star$ and $\omega$,
\begin{align*}
\omega(f, g)h\star k = q(\tau(f)\tau(g)\tau(h)k + \tau(f)\tau(g)\tau(k)h)
\end{align*}
and using this its mediator can be calculated as
\begin{align}
\label{starmediator}
\begin{aligned}
[f, g, h, k] &= q(\omega(f, g)h\star k + \omega(h, k)f\star g - \omega(f, h)g\star k - \omega(g, k)f \star h) = 0,
\end{aligned}
\end{align}
for $f, g, h, k \in \field[q]$.
\end{proof}

\begin{theorem}\label{thm:peircehomomorphism}
The Peirce operator $\rho: \Nring{\field}{\{x\}} \to \field[q]$ is a homomorphism with respect to the commutative product $\star$ defined on $\field[q]$ in Lemma \ref{lem:starmedial}.
\end{theorem}
\begin{proof}
Combine Lemma \ref{peircerecursionlemma} with Corollary \ref{cor:evaluationP1}.
\end{proof}

Although no further use is made here of these observations, Lemma \ref{lem:starfacts} records some interesting facts about the medial algebra $(\field[q], \star)$.
\begin{lemma}\label{lem:starfacts}
\noindent
\begin{enumerate}
\item The medial algebra $(\field[q], \star)$ contains no nonzero idempotent.
\item The subspaces $\balg^{k} = \{f \in \field[q]: \deg(f) \leq k\}$ defined for $k \geq 0$ constitute an ascending filtration of $\field[q]$ satisfying $\balg^{k}\star \balg^{l} \subset \balg^{\max\{k, l\} + 1}$.
\end{enumerate}
\end{lemma}
\begin{proof}
If $f$ is idempotent, then $f = f\star f = 2q\tau(f)f$. If $\tau(f) =0$ this implies $f = 0$, while if $\tau(f) \neq 0$ then comparing degrees in $q$ shows $f = 0$.
By definition of $f\star g$, for nonzero $f,g \in \field[q]$,
\begin{align}\label{degstar}
\deg(f\star g) \leq \max\{\deg(f), \deg(g)\} + 1,
\end{align}
where by convention $\deg(0) = -\infty$. If $f \in \balg^{k}$ and $g \in \balg^{l}$, then \eqref{degstar} shows $\balg^{k}\star \balg^{l} \subset \balg^{\max\{k, l\} + 1}$.
\end{proof}

Next there is described the relation between the Peirce polynomial of $P \in  \Nring{\field}{\{x\}}$ and the first linearization of $P$. For this it is necessary to introduce some terminology.

Given a nonassociative polynomial $P  = P(t_{1}, \dots, t_{m}) \in \Nring{\field}{t_{1}, \dots, t_{m}}$, a commutative $\field$-algebra $\alg$ is said to \emph{satisfy the identity $P$} if $P(a_{1}, \dots, a_{m}) = 0$ for all $a_{1}, \dots, a_{m} \in \alg$. Equivalently, it is also said that $\alg$ \emph{satisfies $P = 0$}.

If $v$ is an indeterminate, the inclusion of sets $\mathsf{T} \to \mathsf{T}\cup \{v\}$ induces an algebra embedding $\Nring{\field}{\mathsf{T}} \to\Nring{\field}{\mathsf{T}\cup\{v\}}$, via which $P \in \Nring{\field}{\mathsf{T}}$ can be viewed as an element of $\Nring{\field}{\mathsf{T}\cup\{v\}}$, also denoted $P$.
Given $x \in \mathsf{T}$, the \emph{(formal) derivative in $x$} is the $\field$-linear map $\nabla_{x}:\Nring{\field}{\mathsf{T}\cup\{v\}} \to\Nring{\field}{\mathsf{T}\cup\{v\}}$ uniquely determined by the requirements that $\nabla_{x}(x)=v$, $\nabla_{x}(y) = 0$ if $y \in \mathsf{T}\setminus\{x\}$, $\nabla_{x}(v) = 0$, and that it satisfy the Leibniz rule
\begin{align}
\nabla_{x}(PQ) = \nabla_{x}(P)Q + P\nabla_{x}(Q)
\end{align}
for $P, Q \in \Nring{\field}{\mathsf{T}\cup\{v\}}$. For example,
\begin{align}\label{nablaex1}
&\nabla_{t_{1}} (t_1^3t_2)=2(t_1(t_1v))t_2 + (vt_1^2)t_2,&
& \nabla_{t_{2}} (t_1^3t_2)=t_{1}^3 v.
\end{align}
Alternatively, for $k \geq 1$ there are $\field$-linear operators $\nabla_{x, k}: \Nring{\field}{\mathsf{T}\cup\{v\}} \to\Nring{\field}{\mathsf{T}\cup\{v\}}$ uniquely determined by the identity
\begin{align}\label{nablaxk}
P(\dots, x + v, \dots) = P(\dots, x, \dots) + \sum_{1 \leq k \leq \weight_{x}(P)}(\nabla_{x, k}P)(\dots, v, \dots)
\end{align}
where $\dots$ indicate elements of $\mathsf{T} \setminus \{x\}$, and $\weight_{v}(\nabla_{x, k}P) = k$. When $k = 1$ there holds $\nabla_{x, 1} = \nabla_{x}$ and the subscript indicating the degree is omitted.

More generally, for an indeterminate $w\notin\mathsf{T}$ write $\nabla_x^{w}$ for the
derivation defined as $\nabla_x$ but with $\nabla_x^{w}(x)=w$ (in this notation $\nabla_x$ is $\nabla_x^{v}$).
For distinct $v,w\notin \mathsf{T}$ the \emph{second linearization} of $P$ in $x$ is
\begin{align}\label{nabla2}
\nabla^2_x P := \nabla_x^{w}\nabla_x^{v}P \in \Nring{\field}{\mathsf{T}\cup\{v,w\}}.
\end{align}
It has degree one in each of $v,w$ and is symmetric in $v,w$ (the derivations $\nabla_x^{v}$,
$\nabla_x^{w}$ commute, since their commutator is a derivation vanishing on generators);
equivalently, $\nabla^2_xP(\dots,v,w)$ is the polarization
$\nabla_{x,2}P(\dots,v+w)-\nabla_{x,2}P(\dots,v)-\nabla_{x,2}P(\dots,w)$ of \eqref{nablaxk}.
The higher linearizations $\nabla^m_x P=\nabla_x^{w_m}\cdots\nabla_x^{w_1}P$ are defined similarly.
For a commutative algebra $\alg$, $(\nabla^2_xP)(c,a,b)$ denotes evaluation at $x\mapsto c$,
$v\mapsto a$, $w\mapsto b$. If $\alg$ satisfies $P$ and $|\field|>\deg P$, then $\alg$ satisfies
$\nabla_xP$ and $\nabla^m_xP$.

Proposition~\ref{pro:Peircedef} relates the Peirce polynomial of an identity $P\in  \Nring{\field}{\{x\}}$ satisfied by a commutative algebra with the evaluation of the first linearization of $P$ at an idempotent. It follows from results of \cite[sec.~4]{Tkachev-universality}, in particular Proposition $4.4$ of that reference. See also \cite[Proposition $27$]{varro2020peirce} which formulates Proposition~\ref{pro:Peircedef} in the multivariate setting.

\begin{proposition}[\cite{Tkachev-universality}]\label{pro:Peircedef}
Given $P \in  \Nring{\field}{\{x\}}$, for any commutative $\field$-algebra $\alg$ and any idempotent $c$ in $\alg$ and any $a \in \alg$ there holds
\begin{equation}\label{PEIRCE}
(\nabla_{x}P)(c, a) =\rho(P;L(c))a,
\end{equation}
where $\rho(P; L(c))$ is the evaluation of the Peirce polynomial $\rho(P; q) \in \field[q]$ on the endomorphism $L(c) \in \End(\alg)$.
\end{proposition}

\begin{proof}
By $\field$-linearity it suffices to prove \eqref{PEIRCE} for $P(x) = x^{\alpha} \in \nass(\{x\})$ a monomial. In this case the claim can be proved by an induction on the degree of $\alpha$ as in the proof of Lemma \ref{peircecompositionlemma} below that uses that the first linearization of $P$ satisfies the Leibniz rule. Details can be found in the proof of \cite[Proposition $4.4$]{Tkachev-universality} or \cite[Proposition $27$]{varro2020peirce}.
\end{proof}

\begin{remark}
The terminology \emph{Peirce polynomial} appears in \cite{Osborn}, but refers implicitly to the polynomial in $L(c)$ determined by the left-hand side of \eqref{PEIRCE}; the abstract Peirce polynomial $\rho(P; q)$ is not formalized and no analogue of the formula \eqref{peircetreeformula} appears. In particular, to make sense of \eqref{PEIRCE} it is necessary to assume the algebra contains a nonzero idempotent, something which is not always the case. The Peirce polynomial and the abstract Peirce spectrum are defined in terms of the identity and do not require such an assumption. \qed
\end{remark}

\begin{definition}\label{def:abstractpeircespectrum}
Given $P \in  \Nring{\field}{\{x\}}$, define the \emph{abstract Peirce spectrum} of $P$, $\sigma_\fie(P)$, to be the set of solutions in $\fie$ of $\rho(P;q)=0$. In particular, if $\rho(P; q) = 0$ then $\sigma_\fie(P) = \fie$.
\end{definition}

\begin{proposition}[{\cite{Tkachev-universality}}]\label{pro:peirceidentity}
If the commutative $\field$-algebra $\alg$ satisfies the identity $P = 0$ with $|\field| > \deg P$ and $c$ is an idempotent in $\alg$, then $L(c)$ satisfies $\rho(P;L(c))=0$. In particular, $\rho(P;q)=0$ for any $q\in \sigma(c)\cap \field$; that is
\begin{equation}\label{sigmain}
\sigma(c)\cap \field \subset \sigma_\fie(P) \quad \text{ for all }c\in \Idm(\alg).
\end{equation}
\end{proposition}
\begin{proof}
That $\rho(P;L(c))=0$ when $\alg$ satisfies $P = 0$ is an immediate consequence of \eqref{PEIRCE}.
\end{proof}

\begin{example}
For $P(x)=x^4$, by Example \ref{principalpowersexample}, $\rho(x^4;q)=2q^3+q^2+q$, and
$$
(\nabla_x P)(x, v)=2((xv)x)x+((x^2)v)x+((x^2)x)v,
$$
so that for a commutative algebra $\alg$ and $c \in \Idm(\alg)$ and $a \in \alg$,
\begin{align*}
(\nabla_x P)(c, a) &=2((ca)c)c+((c^2)a)c+((c^2)c)a\\
&=2c(c(ca))+c(ca)+ca=(2L(c)^3+L(c)^2+L(c))a,
\end{align*}
as follows from Proposition \ref{pro:Peircedef}.
\qed\end{example}

\begin{example}
\label{ex:osborn}
It is instructive to calculate the Peirce polynomial of the Osborn identity \eqref{osborn5}. As $\rho(x;q)=1$, applying \eqref{peircerecursion} repeatedly together with \eqref{peircexn} yields
\begin{align}\label{rhodiverse}
&\rho(x^3x^2;q)=2q^3+3q^2,&
&\rho((x^2)^2;q)=4q^2,&
&\rho((x^2)^2x;q)=4q^3+q,
\end{align}
from which there follows
\begin{equation}\label{onehalfm}
\rho(2x^5-3(x^2)^2x+x^3x^2;q)=4q^4-8q^3+5q^2-q=q(q - 1)(2q-1)^2.
\end{equation}
By Proposition \ref{pro:peirceidentity}, this implies that the Peirce spectrum of any commutative algebra $\alg$ satisfying \eqref{osborn5} is contained in $\{0,\tfrac12,1\}$, a property basic in the study of Jordan and power-associative algebras which both satisfy \eqref{osborn5}.
Note, however, that \eqref{onehalfm} shows that the Peirce polynomial of the Osborn identity has a multiple root $\tfrac12$. It follows from \cite{Tkachev-universality} that, in this situation, the fusion rules need not coincide with those of Jordan algebras.

Here the abstract Peirce spectrum has been defined only to be the algebraic set comprising the roots of the Peirce polynomial, without taking into account their multiplicities. The Peirce polynomial of the Osborn identity illustrates that it is necessary to take into account the extra information given by multiplicities. This could be incorporated by defining the abstract Peirce spectrum as a scheme, but as the issue of multiplicity is not explored further here, this has not been developed.
\qed
\end{example}

\section{Evanescent identities}
\begin{definition}
A nonassociative polynomial $P\in \Nring{\field}{\{x\}}$ is \emph{evanescent} if $\rho(P;q)$ is identically zero.
\end{definition}

The idea of evanescence originates in the study of Peirce polynomials in \cite{Tkachev-universality}, where such polynomials were called \emph{degenerated}. The terminology \emph{evanescent} is taken from \cite{varro2020peirce}.

The goal in what follows is to make a detailed analysis of certain specific evanescent polynomials, without pursuing the general theory or addressing broader natural questions, such as the classification of all evanescent polynomials; see, however, \cite{varro2020peirce} for some results in this direction.

\begin{lemma}\label{evanescentp1lemma}
If $P \in \Nring{\field}{\{x\}}$ is evanescent, then $P(1) = 0$.
\end{lemma}

\begin{proof}
By Corollary \ref{cor:evaluationP1}, $P(1) = \rho(P; \tfrac{1}{2}) = 0$, the last equality by evanescence of $P$.
\end{proof}

\begin{example}
Lemma \ref{evanescentp1lemma} means that the base field $\field$, viewed as a commutative algebra, satisfies the identity $P = 0$ for any homogeneous evanescent $P\in \Nring{\field}{\{x\}}$. This means that any homogeneous evanescent identity is satisfied by some nontrivial commutative algebra.
\qed
\end{example}

\begin{lemma}[{\cite[Proposition $30$]{varro2020peirce}}]\label{lem:evanideal}
The subspace $\evan := \ker \rho \subset \Nring{\field}{\{x\}}$ of evanescent polynomials is an ideal.
\end{lemma}

\begin{proof}
This follows from \eqref{peircerecursion} of Lemma \ref{peircerecursionlemma} and that, by Lemma \ref{evanescentp1lemma}, for $P \in \evan$ there holds $P(1) = 0$.
\end{proof}

Next there are given some examples related with evanescent polynomials.

Applying Definition \ref{defmediator} in the case where the algebra $\alg$ is $ \Nring{\field}{\{x\}}$ yields the \emph{mediator} of nonassociative polynomials $P_{1}, P_{2}, P_{3}, P_{4} \in \Nring{\field}{\{x\}}$:
\begin{align}\label{mediatordefined}
[P_{1}, P_{2}, P_{3}, P_{4}] = (P_{1}P_{2})(P_{3}P_{4}) - (P_{1}P_{3})(P_{2}P_{4}).
\end{align}
\begin{proposition}\label{pro:medial}
For any  $P_{1}, P_{2}, P_{3}, P_{4} \in \Nring{\field}{\{x\}}$, the mediator $[P_{1}, P_{2}, P_{3}, P_{4}]$ is evanescent.
\end{proposition}
\begin{proof}
This follows from Theorem \ref{thm:peircehomomorphism}, which shows that $\rho$ is a homomorphism into the algebra $(\field[q], \star)$ defined in Lemma \ref{lem:starmedial}, and the fact, shown in Lemma \ref{lem:starmedial} that $(\field[q], \star)$ is medial.
\end{proof}

\begin{example}[Medial evanescence]\label{ex:medial}
For any $x^\alpha, x^\beta, x^\gamma, x^\delta \in \Nring{\field}{\{x\}}$,
\begin{equation}\label{medialevan}
\rho((x^\alpha x^\beta)(x^\gamma x^\delta)-(x^\alpha x^\gamma)(x^\beta x^\delta); q)=0,
\end{equation}
follows from the special case of Proposition \ref{pro:medial} with $P_{1}(x) = x^{\alpha}$, $P_{2}(x) = x^{\beta}$, $P_{3}(x) = x^{\gamma}$, and $P_{4}(x) = x^{\delta}$.
\qed\end{example}

\begin{example}
Taking $P_{1}(x)  = x = P_{2}(x)$ and $P_{3}(x) = x^{2} = P_{4}(x)$ in Proposition \ref{pro:medial} yields the evanescence of the identity $E_1$ of \eqref{strange}.
\qed\end{example}

\begin{lemma}
For the \emph{associator} of $P_{1}, P_{2}, P_{3}\in  \Nring{\field}{\{x\}}$ defined by $[P_{1}, P_{3}, P_{2}] := (P_{1}P_{3})P_{2} - P_{1}(P_{3}P_{2})$, there holds
\begin{align}\label{rhoassociator}
\begin{aligned}
\rho([P_{1}, P_{3}, P_{2}]; q) &= (q^{2} - q)P_{3}(1)\,(P_{2}(1)\rho(P_{1};q) -P_{1}(1)\rho(P_{2}; q))\\
&= -(q^{2}-q)P_{3}(1)\Delta(P_{1}, P_{2}; q),
\end{aligned}
\end{align}
where
\begin{align}\label{deltapq}
\Delta(P,Q;t)=\begin{vmatrix} P(1) & Q(1) \\ \rho(P; t) & \rho(Q; t) \end{vmatrix}.
\end{align}
\end{lemma}

\begin{proof}
This follows from Theorem \ref{thm:peircehomomorphism}, which shows that $\rho$ is a homomorphism into the algebra $(\field[q], \star)$ defined in Lemma \ref{lem:starmedial}, and the expression \eqref{starassociator} for the associator of $(\field[q], \star)$.
\end{proof}

\begin{definition}
Define
$\Tri(P, Q, R)$ of $P, Q, R \in \Nring{\field}{\{x\}}^{\times 3}$ by
\begin{align}
\begin{aligned}
\Tri(P, Q, R) &= [P_1, Q_2, R_3] + [P_2, Q_3, R_1] + [P_3, Q_1, R_2]\\
&\quad +  [R_1, Q_2, P_3] + [R_2, Q_3, P_1] + [R_3, Q_1, P_2]\\
& = -P_{1}(Q_{2}R_{3} - Q_{3}R_{2}) - P_{2}(Q_{3}R_{1} - Q_{1}R_{3}) - P_{3}(Q_{1}R_{2} - Q_{2}R_{1})\\
&\quad -R_{1}(Q_{2}P_{3} - Q_{3}P_{2}) - R_{2}(Q_{3}P_{1} - Q_{1}P_{3}) - R_{3}(Q_{1}P_{2} - Q_{2}P_{1}),
\end{aligned}
\end{align}
where $P = (P_{1}, P_{2}, P_{3})$ and similarly for $Q$ and $R$.
\end{definition}

\begin{proposition}\label{pro:associatorevanescence}
For $P, Q, R \in \Nring{\field}{\{x\}}$, $\Tri(P, Q, R)$ satisfies
\begin{align}\label{rhotri}
\begin{aligned}
\rho(&\Tri(P, Q, R) ;q)\\
 &
 = (q^{2} - q) \left( Q_1(1)\Delta(P_2, R_3; q) + Q_2(1)\Delta(P_3, R_1; q) + Q_3(1)\Delta(P_1, R_2; q) \right.\\
&
 \qquad\qquad\quad\left. +Q_1(1)\Delta(R_2, P_3; q) + Q_2(1)\Delta(R_3, P_1; q) + Q_3(1)\Delta(R_1, P_2; q) \right),
\end{aligned}
\end{align}
where $p_{i} = P_{i}(1)$, $q_{i} =Q_{i}(1)$, and $r_{i} = R_{i}(1)$. In particular, $\Tri(P, Q, R)$ is evanescent if the evaluations $p = (p_{1}, p_{2}, p_{3}), q = (q_{1}, q_{2}, q_{3}), r= (r_{1}, r_{2}, r_{3}) \in \field^{3}$ span a subspace of $\field^{3}$ of dimension at most $1$.
\end{proposition}
\begin{proof}
The formula \eqref{rhotri} follows from rewriting
\begin{align}\label{rhotri0}
\begin{aligned}
\rho(&\Tri(P, Q, R);q)\\&= (q^{2} - q)\left((q_{2}p_{3} - q_{3}p_{2})\rho(R_{1}; q) +(q_{3}p_{1} - q_{1}p_{3})\rho(R_{2}; q)+ (q_{1}p_{2} - q_{2}p_{1})\rho(R_{3}; q)  \right)\\
 &  + (q^{2} - q)\left((q_{2}r_{3} - q_{3}r_{2})\rho(P_{1}; q) +(q_{3}r_{1} - q_{1}r_{3})\rho(P_{2}; q)+ (q_{1}r_{2} - q_{2}r_{1})\rho(P_{3}; q)  \right),
\end{aligned}
\end{align}
which is a direct consequence of \eqref{rhoassociator}. The claim about evanescence is evident from \eqref{rhotri0}.
\end{proof}

\begin{example}
Taking $P_{3} = P_{2}$, $Q_{1} = 0$, and $R = P$ in Proposition \ref{pro:associatorevanescence} yields that
\begin{align}
\label{associatordifference1}
&[P_{1}, Q_{2}, P_{2}] + [P_{2}, Q_{3}, P_{1}]
\end{align}
is evanescent provided $Q_{2}(1) = Q_{3}(1)$. Taking $Q_{1} = Q_{2} = Q_{3} = P_{4}$ and $R = P$ in Proposition \ref{pro:associatorevanescence} yields that
\begin{align}
\label{associatordifference2}
&[P_{1}, P_{4}, P_{2}] + [P_{2}, P_{4}, P_{3}] + [P_{3}, P_{4}, P_{1}],
\end{align}
is evanescent, provided $P_{1}(1) = P_{2}(1) = P_{3}(1)$.
\qed\end{example}

A special feature of $\Nring{\field}{\{x\}}$ is that it admits the composition of univariate nonassociative polynomials $P, Q\in \Nring{\field}{\{x\}}$ defined by $(P\circ Q)(x) :=P(Q(x))$.

\begin{lemma}\label{peircecompositionlemma}
Let $Q\in \Nring{\field}{\{x\}}$. If $P\in \Nring{\field}{\{x\}}$ is homogeneous of degree $d\geq 1$, then
\begin{align}\label{PQid}
\rho(P\circ Q; q) = (Q(1))^{d-1}\rho(P; q)\rho(Q; q).
\end{align}
In particular, if $Q(1)=1$ then $\rho(P\circ Q;q)=\rho(P;q)\rho(Q;q)$ for all $P \in \Nring{\field}{\{x\}}$.
\end{lemma}
\begin{proof}
Let $t = Q(1)$. By linearity it suffices to take $P = x^\alpha$ with $|\alpha| = d$. Because evaluation at $1$ is a homomorphism, $Q(x)^\beta(1) = t^{|\beta|}$. The case $d=1$ is trivial. If $x^\alpha = x^\beta x^\gamma$ with $|\beta| = k$, $|\gamma| = m$, then by \eqref{peircerecursion} and the inductive hypothesis,
\begin{align*}
\rho(Q(x)^\alpha;q) &= q\big(t^{m}\rho(Q(x)^\beta;q) + t^{k}\rho(Q(x)^\gamma;q)\big)\\
&= q\big(t^{m}t^{k-1}\rho(x^\beta;q) + t^{k}t^{m-1}\rho(x^\gamma;q)\big)\rho(Q;q)\\
&= t^{d-1}\rho(x^\alpha;q)\rho(Q;q).
\end{align*}
\end{proof}

\begin{proposition}\label{pro:comp}
If $P, Q \in \Nring{\field}{\{x\}}$ satisfy $P(1)=Q(1)=1$, then $P\circ Q - Q\circ P$ is evanescent.
\end{proposition}
\begin{proof}
By Lemma \ref{peircecompositionlemma}, $\rho(P\circ Q; q)= \rho(P;q)\rho(Q;q) = \rho(Q\circ P;q)$.
\end{proof}

\begin{example}[Palintropic evanescence]\label{palintropicexample}
For any nonassociative monomials $x^\alpha, x^\beta\in \Nring{\field}{\{x\}}$,
$$\rho((x^\alpha)^\beta-(x^\beta)^\alpha;q)=0,$$
follows from the special case of Proposition \ref{pro:comp} with $P(x) = x^{\alpha}$ and $Q(x) = x^{\beta}$.
\qed\end{example}

The following property extends \cite[Proposition 8.3]{Tkachev-universality}.

\begin{corollary}\label{pro:inovar}
If $P\in \Nring{\field}{\{x\}}$ is evanescent, then so are $P(x^\alpha)$ and $P(x)^{\alpha}$
for any nonassociative power $\alpha$.
\end{corollary}

\begin{proof}
By Lemma \ref{peircecompositionlemma}, $\rho(P(x^\alpha);q)=\rho(P;q)\rho(x^\alpha;q) = \rho(P(x)^{\alpha};q)$, and all equal $0$ because $\rho(P;q)=0$.
\end{proof}

\section{The second order Peirce operator and fusion rules}
The second-order Peirce operator defined in \cite{Tkachev-universality} arises from evaluating the second linearization of a nonassociative univariate polynomial on elements of a commutative algebra one of which is an idempotent. It can be interpreted as providing fusion rules for eigenvectors of multiplication by this idempotent. Theorem \ref{thm:secondorderpeirce} defines the second-order Peirce operator on $\Nring{\field}{\{x\}}$ via a recursive definition as in the definition in Lemma \ref{peirceoperatorlemma} of the first-order Peirce operator. The analogue for the second-order Peirce operator of the formula \eqref{peircetreeformula} is the formula \eqref{eq:tree_D}, which shows how to calculate the second-order Peirce operator of a nonassociative monomial directly from the associated binary rooted tree.

For vertices $l_1$ and $l_2$ of an abstract binary rooted tree the \emph{nearest common ancestor} $\NCA(l_1, l_2)$ is the root of the minimal subtree of $T$ that contains both $l_1$ and $l_2$.

\begin{theorem}\label{thm:secondorderpeirce}
There is a unique $\field$-linear map $\frak{D}:\Nring{\field}{\{x\}}\to \fie[q, r, s]$ satisfying
\begin{align}
\begin{aligned}\label{defDD}
\mathrm{(i)}\,\, &\frak{D}(x;q, r, s)=0,\\
\mathrm{(ii)}\,\, &\frak{D}(x^{\beta}  x^{\gamma}  ;q, r, s) = q (\frak{D}(x^{\beta}  ;q, r, s)+\frak{D}(x^{\gamma}  ;q,r, s)) \\
&\qquad \qquad \qquad +\rho(x^{\beta};r)\rho(x^{\gamma};s) +\rho(x^{\beta};s)\rho(x^{\gamma};r),
\end{aligned}
\end{align}
for monomials $x^{\beta}, x^{\gamma} \in \Nring{\field}{\{x\}}$. Moreover, for the monomial $x^{\alpha} \in \Nring{\field}{\{x\}}$ associated with the abstract binary rooted tree $T$ having set of leaves $L(T)$, $\frak{D}(x^{\alpha}; q, r, s)$ equals the sum over all unordered pairs of distinct leaves,
\begin{equation}\label{eq:tree_D}
\frak{D}(x^{\alpha}; q, r_{1}, r_{2}) = \sum_{\{l_{1}, l_{2}\} \subseteq L(T)} q^{d(l_{1}, l_{2})} \left( r_{1}^{\height(l_{1}|l_{2}) - 1}r_{2}^{\height(l_{2}|l_{1}) - 1} + r_{2}^{\height(l_{1}|l_{2}) - 1}r_{1}^{\height(l_{2}|l_{1}) - 1} \right),
\end{equation}
where $d(l_{1}, l_{2}) = \height(\NCA(l_{1}, l_{2}))$ denotes the height of the nearest common ancestor $\NCA(l_{1}, l_{2})$ and $\height(l_{1}|l_{2}):= \height(l_1) - d(l_1, l_2)$ is the height of $l_{1}$ relative to $l_2$.
\end{theorem}

Observe that the dependence on the indeterminate $q$ is qualitatively different from the dependence on $r$ or $s$.

\begin{proof}
The map $\frak{D}$ is defined by extending \eqref{eq:tree_D} $\field$-linearly to $\Nring{\field}{\{x\}}$. It suffices to check that $\frak{D}$ so defined satisfies \eqref{defDD}. This is shown by induction on the degree of the binary tree $T$. In the proof it is convenient to write $\frak{D}(x^\alpha)$ in place of $\frak{D}(x^\alpha; q, r, s)$ and to write tree for abstract binary rooted tree.

If $T$ is the trivial tree representing the monomial $x$ then the leaf set contains only $x$, so the set of unordered pairs of distinct leaves is empty and the sum in \eqref{eq:tree_D} and the initial condition $\frak{D}(x) = 0$ is satisfied.

Suppose that $T$ has height at least one and suppose the claim holds for all subtrees of degree lower than that of $T$. Then $T$ can be written in a unique way as the graft, $T = A \vee B$, of two subtrees $A$ and $B$ each of degree lower than $T$. Denote by $\binom{L(T)}{2}$ the set of unordered pairs of leaves of $T$. The set of leaves of $T$ is the disjoint union of the sets of leaves of its subtrees: $L(T) = L(A) \cup L(B)$. Consequently,
\begin{align}\label{leafpairs}
\binom{L(T)}{2} = \binom{L(A)}{2} \cup \binom{L(B)}{2} \cup \Big\{ \{l,r\} \mid l \in L(A), r \in L(B) \Big\}.
\end{align}
Let $x^\alpha$ and $x^\beta$ be the monomials corresponding with $A$ and $B$ so that the monomial corresponding with $T$ is $x^{\alpha}x^{\beta}$. The inductive hypothesis implies that \eqref{eq:tree_D} holds for each of $x^\alpha$ and $x^\beta$. It will be used to calculate the right-hand side of the recursive formula \eqref{defDD} and this will be shown to equal \eqref{eq:tree_D} for $T$.

(i) Consider pairs of leaves contained in the subtree $A$, $\{l_1, l_2\} \subseteq L(A)$. The root of $A$ has height $1$ relative to the root of $T$, so the height of the nearest common ancestor of $l_1, l_2$ relative to the root of $T$ is one more than its height relative to the root of $A$, $d_T(l_1, l_2) = d_A(l_1, l_2) + 1$. Similarly, the height in $T$ of a leaf of $A$ is one more than its height in $A$, $\height_T(v) = \height_A(v) + 1$. Consequently,
\begin{align}\label{subtreeinvariance}
\begin{aligned}
\height_T(l_1|l_2) -1 &= \height_T(l_1) - d_T(l_1, l_2) - 1 = (\height_A(l_1) + 1) - (d_A(l_1, l_2) + 1) - 1 \\
&= \height_A(l_1) - d_A(l_1, l_2) - 1 = \height_A(l_1|l_2) -1,
\end{aligned}
\end{align}
so the relative height of $l_1$ with respect to $l_2$ is the same in $A$ as in $T$, and the same is true with $l_2$ and $l_1$ interchanged.
Therefore, the relative path distances from the local ancestor to the leaves remain invariant.
Multiplying by $q$ yields $q^{d_A(l_1, l_2)}$ to $q^{d_A(l_1, l_2)+1} = q^{d_T(l_1, l_2)}$.
By the inductive hypothesis and \eqref{subtreeinvariance}
\begin{align}\label{partialA}
\begin{aligned}
q\frak{D}&(x^{\alpha}; q, r_{1}, r_{2})\\
&= q\sum_{\{l_{1}, l_{2}\} \subseteq L(A)} q^{d_A(l_{1}, l_{2})} \big( r_{1}^{\height_A(l_{1}|l_{2}) - 1}r_{2}^{\height_A(l_{2}|l_{1}) - 1} + r_{2}^{\height_A(l_{1}|l_{2}) - 1}r_{1}^{\height_A(l_{2}|l_{1}) - 1} \big)\\
&= \sum_{\{l_{1}, l_{2}\} \subseteq L(T)} q^{d_T(l_{1}, l_{2})} \left( r_{1}^{\height_T(l_{1}|l_{2}) - 1}r_{2}^{\height_T(l_{2}|l_{1}) - 1} + r_{2}^{\height_T(l_{1}|l_{2}) - 1}r_{1}^{\height_T(l_{2}|l_{1}) - 1} \right).
\end{aligned}
\end{align}
The same argument establishes the parallel formula with $B$ and $\beta$ in place of $A$ and $\alpha$.

(ii) Consider pairs $\{l_1, l_2\}$ such that $l_1 \in L(A)$ and $l_2 \in L(B)$. Because these leaves belong to separate subtrees, their nearest common ancestor is the root of $T$, so that $d_T(l_1, l_2) = 0$ and $q^{d_T(l_1, l_2)} = 1$. The heights of $l_1$ and $l_2$ in $A$ and $B$ are one less than their heights in $T$, so that $\height_{A}(l_1) = \height_T(l_1) - 1$ and $\height_B(l_2) = \height_T(l_2) - 1$. Consequently $\height_{T}(l_1|l_2) = \height_T(l_1) - d_T(l_1, l_2) = \height_A(l_1) + 1$ and similarly $\height_{T}(l_2|l_1) =  \height_B(l_2) + 1$.
By the formula \eqref{peircetreeformula} for the Peirce operator $\rho$,
 \begin{align*}
 \rho(A;r)\rho(B;s) &= \Big( \sum_{l_1 \in L(A)} r^{\height_A(l_1)} \Big) \Big( \sum_{l_2 \in L(B)} s^{\height_B(l_2)} \Big) \\
 &= \sum_{l_1 \in L(A), l_2\in L(B)} r^{\height_T(l_1|l_2)-1}s^{\height_T(l_2|l_1)-1}\,.
\end{align*}
Symmetrizing this in $r$ and $s$ yields
\begin{align}\label{rhoarhob}
\begin{aligned}
\rho(A;r)&\rho(B;s) + \rho(A;s)\rho(B;r) \\
&= \sum_{l_1 \in L(A), l_2 \in L(B)}q^{d_T(l_1, l_2)} \Big( r^{\height_T(l_1|l_2)-1}s^{\height_T(l_2|l_1)-1} + s^{\height_T(l_1|l_2)-1}r^{\height_T(l_2|l_1)-1}\Big).
\end{aligned}
\end{align}
In conjunction with \eqref{leafpairs}, summing \eqref{partialA} and its parallel with $B$ in place of $A$ with \eqref{rhoarhob} yields that $q (\frak{D}(x^{\alpha}  ;q, r, s)+\frak{D}(x^{\beta}  ;q,r, s)) +\rho(x^{\alpha};r)\rho(x^{\beta};s) +\rho(x^{\alpha};s)\rho(x^{\beta};r)$ equals the sum on the right-hand side of \eqref{eq:tree_D}. This completes the inductive step and (after relabeling) shows the claim.
\end{proof}

\begin{definition}\label{def:secondorderpeirce}
$\frak{D}(P;q, r, s)$ is the \emph{second order Peirce operator} of $P \in \Nring{\field}{\{x\}}$.
\end{definition}

\begin{corollary}
$\deg_q \frak{D}(P;q, r, s)\le \deg_{x}(P)-2$\,.
\end{corollary}
\begin{proof}
This follows from \eqref{eq:tree_D} because the height of the nearest common ancestor of two leaves in an abstract binary rooted tree is at most $2$ less than the number of leaves of the tree (which equals the degree of the corresponding monomial).
\end{proof}

\begin{lemma}\label{secondpeircerecursionlemma}
For all $P, Q \in \Nring{\field}{\{x\}}$ there holds
\begin{align}
\begin{aligned}
\label{secondpeircerecursion}
\frak{D}(P  Q  ;q, r, s)
&=q(Q(1)\frak{D}(P ;q, r, s)+P(1)\frak{D}(Q ;q, r, s))\\
&+\rho(P;r)\rho(Q;s) +\rho(P;s)\rho(Q;r)
\end{aligned}
\end{align}
\end{lemma}
\begin{proof}
This follows from \eqref{defDD} by $\fie$-linearity of the second order Peirce operator.
\end{proof}

\begin{example}[Principal and plenary powers]
Let $T$ be the  binary tree representing the fourth principal power $x^4=x(x(xx))$ (see the left part of Fig.~\ref{fig:forest_x4}).
\begin{figure}[ht]
\begin{center}
\begin{tikzpicture}[
  level distance=1cm,
  level 1/.style={sibling distance=13mm},
  level 2/.style={sibling distance=13mm},
  level 3/.style={sibling distance=13mm}
    ]
  \node[style={circle}]
  {\XX}
    child {node {$l_1$}}
    child {node {$\bullet$}
    child {node {$l_2$} }
    child {node {$\bullet$} child {node {$l_3$}} child {node {$l_4$}}}
    }

;
\end{tikzpicture}
\qquad \qquad
\begin{tikzpicture}[
  level distance=1cm,
  level 1/.style={sibling distance=17mm},
  level 2/.style={sibling distance=9mm},
  level 3/.style={sibling distance=7mm}
    ]
  \node[style={circle}]
{\XX}
    child {\NN
    child {node {$m_1$}}
    child {node {$m_2$} }
    }
    child {\NN
    child {node {$m_3$} }
    child {node {$m_4$} }
    }

;
\end{tikzpicture}
\end{center}
\caption{Binary trees representing principal $x^4$ and plenary $x^2x^2$ powers}
    \label{fig:forest_x4}
\end{figure}
The tree has four leaves labeled $l_1, l_2, l_3, l_4$ having heights $\height(l_1) = 1$, $\height(l_2) = 2$, $\height(l_3) = \height(l_4) = 3$ and relative heights
\begin{align*}
&\height(l_1|l_2) = 1 -0 = 1, &&\height(l_1|l_3) = 1-0 = 1, &&\height(l_1|l_4) = 1-0 = 1,&\\
&\height(l_2|l_3) = 2-1 = 1,&&\height(l_2|l_4) = 2-1 = 1,&&\height(l_3|l_4) = 3 -2 = 1,\\
&\height(l_2|l_1) = 2-0 = 2, &&\height(l_3|l_1) = 3-0 = 3, &&\height(l_4|l_1) = 3-0 = 3,&\\
&\height(l_3|l_2) = 3- 1= 2,&&\height(l_4|l_2) = 3 - 1 =2,&&\height(l_4|l_3) = 3 -2 = 1.
\end{align*}
Evaluating the $\binom{4}{2} = 6$ unordered leaf pairs according to the depth of the nearest common ancestor $d(l_i, l_j)$ yields three cases:
\begin{enumerate}
    \item The pair $\{l_3, l_4\}$ satisfying $d(l_3, l_4) = 2$ contributes
    $ q^2 \left(r^{1-1}s^{1-1} + s^{1-1}r^{1-1}\right) = 2q^2$.

    \item Each of the pairs $\{l_2, l_3\}$ and $\{l_2, l_4\}$ with $d(l_i, l_j) = 1$ contributes
    \begin{align*}
       q^1 \left(r^{1-1}s^{2-1} + s^{1-1}r^{2-1}\right) = q(r+s).
    \end{align*}
    Combining these two pairs gives $2q(r+s)$.

    \item Each of the pairs containing the leaf $l_1$, namely $\{l_1, l_2\}$, $\{l_1, l_3\}$, and $\{l_1, l_4\}$, has $d(l_1, l_i) = 0$ and their contributions are:
    \begin{align*}
        \{l_1, l_3\}, \, \{l_1, l_4\} &\implies q^0 \left(r^{3-1}s^{1-1} + s^{3-1}r^{1-1}\right)  = r^2 + s^2 \\
                \{l_1, l_2\} &\implies q^0 \left(r^{2-1}s^{1-1} + s^{2-1}r^{1-1}\right)  = r + s
    \end{align*}
    Summing these gives $2(r^2 + s^2) + (r + s)$.
    \end{enumerate}
Summing the contributions of the $6$ pairs yields
\begin{equation}\label{dfourthpower}
    \frak{D}(x^4; q, r, s) = 2q^2 + 2q(r+s) + 2(r^2 + s^2) + (r+s).
\end{equation}
Arguing similarly for the plenary power yields
\begin{equation}\label{dplenarypower}
\frak{D}(x^2x^2; q, r, s) = 4q + 8rs.
\end{equation}
\qed\end{example}

Proposition \ref{pro:Peir2nd} is for the second-order Peirce operator what Proposition \ref{pro:Peircedef} is for the first-order Peirce operator.

\begin{proposition}[{\cite[Proposition~5.1]{Tkachev-universality}}]\label{pro:Peir2nd}
Given $P \in  \Nring{\field}{\{x\}}$, for any commutative $\field$-algebra $\alg$ and any idempotent $c$ in $\alg$ and any $a\in \alg_c(\lambda)$ and $b\in \alg_c(\mu)$ there holds
\begin{equation}\label{PEIRCE1}
(\nabla^2_{x}P)(c, a,b) =\frak{D}(P;L(c),\lambda,\mu)(ab).
\end{equation}
where $\frak{D}(P;L(c),\lambda,\mu)$ is the evaluation of the second-order Peirce operator $\frak{D}(P;q, r, s)\in \field[q, r, s]$ sending $q$ to the endomorphism $L(c) \in \End(\alg)$ and $r$ and $s$ to the scalars $\lambda, \mu \in \field$.
\end{proposition}

\begin{proof}
This follows from \cite[Proposition~5.1]{Tkachev-universality}.
\end{proof}

The relation \eqref{PEIRCE1} contains information about the fusion rules associated with the identity $P=0$. For a commutative algebra satisfying $P=0$, the left-hand side of \eqref{PEIRCE1} vanishes, implying that $\frak{D}(P;L(c),\lambda,\mu)$ vanishes on $\alg_c(\lambda)\alg_c(\mu)$. This observation underlies Lemma~\ref{lem:2ndPeirce}.

\begin{definition}\label{def:abstractfusionrules}
Given $P \in  \Nring{\field}{\{x\}}$, define the \emph{abstract fusion rules} of $P$, $\fusion_\fie(P)$, to be the algebraic set
\begin{align}
\{(\alpha, \lambda, \mu) \in \fie^{3}: \alpha, \lambda, \mu \in \sigma_{\field}(P), \frak{D}(P;\alpha, \lambda, \mu)=0\}.
\end{align}
In particular, if $\frak{D}(P; q, r, s) = 0$ then $\fusion_\fie(P) = \sigma_{\field}(P)^{3}$.
\end{definition}

Lemma~\ref{lem:2ndPeirce} shows that if the commutative $\field$-algebra $\alg$ satisfies the identity $P$ then any set $(\alpha, \lambda, \mu)$ for which $\alg_{c}(\lambda)\alg_{c}(\mu)$ is contained in a sum of eigenspaces including $\alg_{c}(\alpha)$ is contained in $\fusion_{\field}(P)$. In this sense the concrete fusion rules satisfied by idempotents in $\alg$ are described in terms of parameters contained in $\fusion_{\field}(P)$, so that the abstract fusion rules of $P$ provide a priori constraints on the fusion rules satisfied by any commutative algebra satisfying $P$.

\begin{lemma}\label{lem:2ndPeirce}
Suppose $|\field| > \deg P$. If the commutative $\field$-algebra $\alg$ satisfies the identity $P = 0$, $c$ is an idempotent in $\alg$, and $\lambda,\mu\in \sigma_\fie(P)$, then
\begin{equation}\label{fusionrules0}
\alg_c(\lambda)\alg_c(\mu)\subset \ker \frak{D}(P;L(c),\lambda,\mu).
\end{equation}
If, moreover, $\frak{D}(P;q,\lambda,\mu)$ is a nonzero polynomial in $q$ that is a product of distinct linear factors over $\field$, or if $L(c)$ is diagonalizable over $\field$, then there holds the fusion rule
\begin{align}\label{fusionrules}
\begin{aligned}
\alg_c(\lambda)\alg_c(\mu)&\subset \bigoplus_{\alpha\in \Omega_P(\lambda,\mu)}\alg_c(\alpha),\\
 & \text{where }\,\,\Omega_P(\lambda,\mu)=\{\alpha \in \sigma_{\fie}(P): \frak{D}(P;\alpha,\lambda,\mu)=0\}.
\end{aligned}
\end{align}
There are two mutually exclusive possibilities for the index set $\Omega_P(\lambda,\mu)$:
\begin{enumerate}
[label=$\mathrm{(D\arabic*)}$,itemsep=0.5ex,leftmargin=1cm]
\item
If $\frak{D}(P;\alpha,\lambda,\mu)=0$ then $\Omega_P(\lambda,\mu)=\sigma_{\fie}(P)$.
\item
If $\frak{D}(P;\alpha,\lambda,\mu)$ is nonzero and has degree $d\ge 0$ in $\alpha$ then the cardinality of $\Omega_P(\lambda,\mu)$ is at most $d$ (if $d = 0$ this means that $\alg_c(\lambda)\alg_c(\mu) = \{0\}$).
\end{enumerate}
\end{lemma}

\begin{proof}
Because $|\field| > \deg P$, the algebra $\alg$ satisfies $\nabla^2_xP$, so the left-hand side of \eqref{PEIRCE1} vanishes for $a \in \alg_c(\lambda)$ and $b \in \alg_c(\mu)$, and Proposition \ref{pro:Peir2nd} yields \eqref{fusionrules0}; see also \cite[Theorem $6.4$]{Osborn}. If $\alpha \in \field$ and $\ker(L(c) - \alpha\Id) \neq \{0\}$, then $\alpha$ is an eigenvalue of $L(c)$, so $\rho(P;\alpha) = 0$ by Proposition \ref{pro:peirceidentity}, that is, $\alpha \in \sigma_\fie(P)$.
If $\frak{D}(P;q,\lambda,\mu)=\gamma\prod_i(q-\alpha_i)$ with $\gamma \neq 0$ and the $\alpha_i\in\field$ pairwise distinct, the factors $q - \alpha_i$ are pairwise coprime, so $\ker\frak{D}(P;L(c),\lambda,\mu)=\bigoplus_i\ker(L(c)-\alpha_i\Id)$; by the preceding observation, only $\alpha_i \in \Omega_P(\lambda,\mu)$ contribute nonzero summands, which gives \eqref{fusionrules}.
If $\alg = \bigoplus_{\alpha \in \field}\alg_c(\alpha)$, then $\frak{D}(P;L(c),\lambda,\mu)$ acts on $\alg_c(\alpha)$ as multiplication by $\frak{D}(P;\alpha,\lambda,\mu)$, so its kernel is the sum of those $\alg_c(\alpha)$ with $\frak{D}(P;\alpha,\lambda,\mu) = 0$, and again only $\alpha \in \sigma_\fie(P)$ contribute nonzero summands, which gives \eqref{fusionrules}.
The claims \textup{(D1)} and \textup{(D2)} follow from the definition of $\Omega_P(\lambda,\mu)$, and the case $d = 0$ of \textup{(D2)} from \eqref{fusionrules0}, because a nonzero constant has trivial kernel.
\end{proof}

\begin{example}
Consider $P=x^4-(x^2)^2$. By \eqref{peircexn} and \eqref{rhodiverse}, $\rho(P;q)=2q^3-3q^2+q=(2q-1)(q-1)q$, so the Peirce spectrum of $P$ is $\sigma(P)=\{0,\tfrac12,1\}$.
Combining \eqref{dplenarypower} and \eqref{dfourthpower} yields
\begin{align}
\frak{D}(P;q, r, s)&=2q^2 - 4q + (2q+1)(r+s) + 2r^2 - 8rs + 2s^2.
\end{align}
Lemma~\ref{lem:2ndPeirce} yields fusion rules for $(\lambda,\mu)\in \sigma(P)\times \sigma(P)$. For example, for $(\lambda,\mu)=(\tfrac12,\tfrac12)$ there results
$\frak{D}(P;q,\tfrac12,\tfrac12)=2q(q-1)$ so that $\Omega(\tfrac12,\tfrac12)=\{0,1\}$, which implies the well-known fusion rule for power associative algebras: $\alg_c(\tfrac12)\alg_c(\tfrac12)\subset \alg_c(0)\oplus \alg_c(1)$.
\qed\end{example}

\smallskip
\begin{lemma}[Proposition~5.2, \cite{Tkachev-universality}]
\label{lem:Devan}
For any $P \in  \Nring{\field}{\{x\}}$,
\begin{equation}\label{Devan}
\frak{D}(P;q, r,\tfrac12)=\frac{\rho(P;q)-\rho(P;r)}{q-r}.
\end{equation}
\end{lemma}
\begin{proof}
This follows straightforwardly from an induction using \eqref{defDD}.
\end{proof}

The second order Peirce operator applied to an evanescent identity implies fusion rules \eqref{fusionrules} for Peirce eigenspaces of eigenvalues not equal to $\tfrac{1}{2}$. On the other hand, Corollary \ref{cor:D12} shows that the second order Peirce operator places no constraints on fusion rules involving the eigenvalue $\tfrac12$.

\begin{corollary}\label{cor:D12}
If $P \in \Nring{\field}{\{x\}}$ is evanescent, then $\frak{D}(P;q, r,\tfrac12)=0$.
\end{corollary}

\begin{proof}
Immediate from \eqref{Devan}.
\end{proof}

\begin{theorem}
For $P_{1}, P_{2}, P_{3}, P_{4} \in \Nring{\field}{\{x\}}$ and $\alpha, \lambda, \mu \in \field$ the mediator satisfies
\begin{align}\label{DPPPP}
\begin{aligned}
\mathfrak{D}&\big([P_1, P_2, P_3, P_4]; q, r,s) \\
&= (q-rs)\left(\Delta(P_1,P_4;s)\Delta(P_2,P_3;r)+\Delta(P_1,P_4;r)\Delta(P_2,P_3;s)\right),
\end{aligned}
\end{align}
where $\Delta(P,Q;t)$ is defined in \eqref{deltapq}.
\end{theorem}

\begin{proof}
Write $M=[P_1, P_2, P_3, P_4]=M_1-M_2$, where $M_1 = (P_1 P_2) (P_3  P_4)$ and $M_2 = (P_1  P_3)  (P_2  P_4)$. Then  $\mathfrak{D}(M) = \mathfrak{D}(M_1) - \mathfrak{D}(M_2)$  (the indeterminate arguments $q, r, s$ of $\mathfrak{D}$ are omitted for readability). Let $X = P_1 P_2$ and $Y = P_3 P_4$. Then \eqref{secondpeircerecursion} yields
\begin{equation*}
\mathfrak{D}(M_1) = q\big(Y(1)\mathfrak{D}(X) + X(1)\mathfrak{D}(Y)\big) + \rho(X; r)\rho(Y; s) + \rho(X; s)\rho(Y; r)
\end{equation*}
Evaluating at $1$ yields $X(1) = P_1(1)P_2(1)$ and $Y(1) = P_3(1)P_4(1)$. Expanding the inner components recursively yields
\begin{align*}
\mathfrak{D}(X) &= q\big(P_2(1)\mathfrak{D}(P_1) + P_1(1)\mathfrak{D}(P_2)\big) + \rho(P_1; r)\rho(P_2; s) + \rho(P_1; s)\rho(P_2; r) \\
\mathfrak{D}(Y) &= q\big(P_4(1)\mathfrak{D}(P_3) + P_3(1)\mathfrak{D}(P_4)\big) + \rho(P_3; r)\rho(P_4; s) + \rho(P_3; s)\rho(P_4; r)
\end{align*}
The cross-terms can be reorganized to yield
\begin{align*}
\rho(X; r)\rho(Y; s) &= rs \big(P_2(1)\rho(P_1;r) + P_1(1)\rho(P_2;r)\big)\big(P_4(1)\rho(P_3;s) + P_3(1)\rho(P_4;s)\big) \\
\rho(X; s)\rho(Y; r) &= rs \big(P_2(1)\rho(P_1;s) + P_1(1)\rho(P_2;s)\big)\big(P_4(1)\rho(P_3;r) + P_3(1)\rho(P_4;r)\big).
\end{align*}
Arguing similarly for $M_2$ and grouping the remaining nonvanishing terms multiplying $q$ and $rs$ yields
\begin{align}
\begin{aligned}
\mathfrak{D}&\big([P_1, P_2, P_3, P_4]; q, r, s)\\
&= (q-rs)P_1(1)P_2(1) \Big( \rho(P_3; r)\rho(P_4; s) + \rho(P_3; s)\rho(P_4; r) \Big) \\
&+ (q-rs)P_3(1)P_4(1) \Big( \rho(P_1; r)\rho(P_2; s) + \rho(P_1; s)\rho(P_2; r) \Big) \\
&- (q-rs)P_1(1)P_3(1) \Big( \rho(P_2; r)\rho(P_4; s) + \rho(P_2; s)\rho(P_4; r) \Big) \\
&- (q-rs)P_2(1)P_4(1) \Big( \rho(P_1; r)\rho(P_3; s) + \rho(P_1; s)\rho(P_3; r) \Big),
\end{aligned}
\end{align}
and rewriting this yields \eqref{DPPPP}.
\end{proof}

\begin{corollary}\label{cor:Devan}
Suppose $P_{1}, P_{2}, P_{3}, P_{4} \in \Nring{\field}{\{x\}}$, where $P_1$ is an evanescent polynomial. Then $\rho([P_{1}, P_{2}, P_{3}, P_{4}]; q)=0$ and $\mathfrak{D}([P_{1}, P_{2}, P_{3}, P_{4}]; q, r, s)=0$.
\end{corollary}

\begin{proof}
If $P_{1}$ is evanescent, then $[P_{1}, P_{2}, P_{3}, P_{4}]$ is evanescent because evanescent polynomials are an ideal, and by the theorem, $\mathfrak{D}\big([P_1, P_2, P_3, P_4]; q, r, s) = 0$.
\end{proof}

\section{Degree $6$ evanescent polynomials}\label{sec:degreesix}
In this section and those following it is assumed everywhere without further comment that $|\field| > 6$. This guarantees, by the discussion preceding Proposition \ref{pro:Peircedef}, that an algebra satisfying a degree $6$ identity also satisfies its linearizations, so that Lemma \ref{lem:2ndPeirce} applies to $E_1$ and $E_2$.

Let $\evan_n$ denote the vector space of evanescent univariate homogeneous degree $n$ polynomials $P\in \Nring{\field}{\{x\}}$.

Observe that, because the only nonassociative monomial whose Peirce polynomial contains a multiple of $q^{n-1}$ is $x^{n}$, no element of $\evan_{n}$ contains a summand proportional to the principal power $x^{n}$.

Varro \cite[Prop. 43]{varro2020peirce} proved that $\evan_n=\{0\}$ for $n\le5$ and $\dim \evan_n\ge W_n-n+2$ for $n \geq 3$, where $W_n$ is the Wedderburn-Etherington number defined in \eqref{wedder}. Varro proved that the minimum degree $n \geq 3$ for which $W_{n} - n +2 > 0$ is $n = 6$, and for this degree he showed that $\evan_6=\operatorname{span}(E_1, E_2)$, where $E_1$ is defined in \eqref{strange} and
\begin{align}
E_2&=(x^{3})^2+ x(x(x^2)^2)-x^2x^4 - x(x^2x^3),\label{Sprimeekv}
\end{align}
although he worked instead with the basis $\{E_2,E_1 + E_2\}$.

That $E_1$ is evanescent is a special case of Example \ref{palintropicexample}. That $E_{2}$ is evanescent, $\rho(E_1 - E_2; q) = 0$, follows from \eqref{rhoassociator} applied to
\begin{align}\label{SSprimedef}
\begin{aligned}
E_1 - E_2& = (x^2)^3-2(x^{3})^2 - x(x(x^2)^2) + x^2x^4 +x(x^2x^3) \\
&= [x,x, (x^2)^2] + [x^3, x^2, x] + [x^3, x, x^2].
\end{aligned}
\end{align}
It is convenient to give here an elementary proof that $E_1$ and $E_2$ span $\evan_6$.

\begin{proposition}[{\cite[Prop. 43]{varro2020peirce}}]\label{pro:adiabatic}
If $\chr\field \neq 2$, then $\evan_6 = \field\{E_1, E_2\}$.
\end{proposition}

\begin{proof}
The most general nonassociative polynomial of degree $6$ has the form
\begin{align}
P = Ax^{6} + B x^2x^4 + C (x^3)^2 + D(x^2)^3 + Ex(x(x^2)^2) + Fx(x^2x^3).
\end{align}
Applying the Peirce operator yields
\begin{align}
\begin{aligned}
\rho(P; q) &= 2Aq^5 + (A + 2B + 4E + 2F)q^4 + (A + B + 4C + 4D + 3F)q^3\\
&+ (A + 3B + 2C + 2D + E)q^2 + (A + E + F)q.
\end{aligned}
\end{align}
If $\chr \field \neq 2$, this vanishes if and only if $A = 0$, $F = B = -E$, and $B + C + D = 0$. Taking $B = 0$ and $D = -C = 1$ yields $E_1$, while taking $B = -1$, $C = 1$, $E = 1$, and $F = -1$ yields $E_2$.
\end{proof}

Proposition~\ref{pro:adiabatic} implies that no restriction is imposed a priori on the Peirce spectrum of an algebra satisfying a linear combination of $E_1$ and $E_2$.
As an illustration, Example \ref{ex:10} shows that any $\lambda\in\field$ can occur in the Peirce spectrum of an idempotent in a $(2, 3)$-palintropic algebra.

\begin{example}\label{ex:10}
Let $\alg$ be the commutative algebra generated by two elements $e$ and $u$ subject to the relations: $e^2=e$, $eu=\lambda u$ and $u^2=0$, where $\lambda\in \field$. This algebra is not associative if and only if $\lambda \notin\{0, 1\}$. Then for any $x=ae+bu$,
\begin{align*}
&x^2=a^2e+2ab\lambda u,&
&(x^2)^2=a^4e+4a^3b\lambda^2 u,&
(x^2)^3 =a^6e+2a^5b\lambda^2(2\lambda+1) u,
\end{align*}
and similarly
\begin{align*}
&x^3=a^3e+a^2b\lambda(2\lambda+1) u,
&(x^3)^2=a^6e+2a^5b\lambda^2(2\lambda+1) u,
\end{align*}
from which it follows that $\alg$ satisfies $E_1$. Note that $\field\{u\}$ is an ideal of $\alg$ and a trivial subalgebra. By construction, $\lambda$ is in the Peirce spectrum of the idempotent $e$.
\qed\end{example}

The structure of the Peirce spectrum of an algebra satisfying an evanescent identity can be constrained by fusion rules derived from the second linearization of the identity at an idempotent. For the evanescent identity $E_1$ the fusion rules have a particularly simple form.

The identity $E_{1}$ is distinguished within $\evan_6$. The second linearizations of both $E_1$ and $E_2$ at an idempotent $c$  are nontrivial, but behave differently. Proposition~\ref{pro:de1e2} formalizes the observation that the second linearization of $E_1$ is linear in $L(c)$, whereas the second linearization of $E_2$ is quadratic in $L(c)$. In conjunction with Lemma~\ref{lem:2ndPeirce} this yields Proposition~\ref{pro:fusion} showing that algebras satisfying $E_{1}$ have very simple fusion rules, which are studied further in Section \ref{sec:adiabatic}.

\begin{proposition}\label{pro:de1e2}
For $E_{1}, E_{2} \in \Nring{\field}{\{x\}}$ as in \eqref{strange} and \eqref{Sprimeekv},
\begin{align}\label{DE1E2}
\begin{aligned}
\frak{D}(E_1;q, r , s )&= 2(2r  - 1)(2s  - 1)(q - r s ), \\
\frak{D}(E_2;q, r , s )&= (2r  - 1)(2s  - 1)\Big[ (q - 1)(2q - r  - s ) - (r  - s )^2 \Big]\\
& = (2r  - 1)(2s  - 1)\Big[2q^2 - (r  + s  + 2)q - (r  - s )^{2} + r  + s \Big].
\end{aligned}
\end{align}
In particular, for any $E \in \evan_6$, $\frak{D}(E;q, r , s )$ is at most quadratic in $q$ and is linear in $q$ if and only if $E$ is a nonzero multiple of $E_1$.
\end{proposition}

\begin{proof}
Using \eqref{defDD} yields
\begin{align}
\begin{aligned}
\frak{D}((x^2)^3;q, r , s )& = 2q(2q + 1) + 8r s (q + r  + s ), \\
\frak{D}((x^{3})^2;q, r , s )&= 4q(q + r  + s ) + 2r s (2r  + 1)(2s  + 1), \\
\frak{D}(x(x(x^2)^2);q, r , s )&= 4q^3 + 8q^2r s  + 4q(r ^2 + s ^2) + 4(r ^3 + s ^3) + (r  + s ) ,\\
\frak{D}(x^2x^4;q, r , s )&= 2q^3 + 2q^2(r  + s ) + q\big(2(r ^2 + s ^2)+ r  + s  + 2\big) \\
&\quad+ 4r s (r ^2 + s ^2) + 2r s (r  + s ) + 4r s , \\
\frak{D}(x(x^2x^3);q, r , s )&= 2q^3 + 2q^2(r  + s  + 1) + 4qr s (r  + s  + 1) \\
&\quad+ 2(r ^3 + s ^3) + 3(r ^2 + s ^2),
\end{aligned}
\end{align}
which when assembled yield \eqref{DE1E2}.
\end{proof}

\begin{proposition}\label{pro:fusion}
If $\alg$ is a commutative $\field$-algebra satisfying $E_{1}$, then in the fusion rules for eigenvalues $\lambda, \mu$ distinct from $\tfrac{1}{2}$ of a nonzero idempotent $c \in \alg$ there is at most a single summand in \eqref{fusionrules}.
\end{proposition}
\begin{proof}
Combine Proposition \ref{pro:de1e2} with \eqref{fusionrules}.
\end{proof}

For a commutative $\field$-algebra satisfying $E \in \evan_6$ not proportional to $E_{1}$ it follows from Proposition \ref{pro:de1e2} that $\frak{D}(E;q, r,s)$ is quadratic in $q$ and this means that for eigenvalues $\lambda, \mu$ distinct from $\tfrac{1}{2}$ of a nonzero idempotent $c \in \alg$ there are at most two summands in the fusion rule \eqref{fusionrules} of Lemma~\ref{lem:2ndPeirce}. In the cases where the discriminant of $\frak{D}(E;q, \lambda, \mu)$ in $q$ vanishes there can in principle appear a single summand in \eqref{fusionrules}, but in general the fusion rules will still be more complicated than in the case of algebras satisfying $E_1$. A complete discussion is omitted because it requires careful consideration of arithmetic properties depending on $\field$ which are not the focus here. In either case for particular algebras some or all of the products $\alg_{c}(\lambda)\alg_{c}(\mu)$ can be trivial. In this sense, $E_1$ is distinguished within $\evan_6$ in that its abstract fusion rules are simpler than those of $E \in \evan_6$ not proportional to $E_1$. Moreover, as is described in detail in Section \ref{sec:adiabatic}, they can be made completely explicit and have an interesting multiplicative structure.

\begin{remark}\label{remark:Osborn}
Osborn observed that the bivariate polynomial
\begin{align}
O(x, y):= 2(y^2x)x - 2((yx)y)x  - 2((yx)x)y + 2(x^2y)y - x^2y^2+ (xy)^2
\end{align}
in \cite[$(5.4)$]{Osborn} has the property that its Peirce polynomial vanishes; substituting $x^{2}$ for $y$ in $O(x, y)$ yields $O(x, x^{2})= E_1(x) + 2E_2(x)$. The set of commutative algebras satisfying $E_1 + 2E_2$ contains properly that satisfying $O$. Osborn calculated the Peirce decomposition associated with $O$ and gave as examples two of its fusion rules \cite[p. 220]{Osborn}. The second-order Peirce polynomial of $O$ has the form as in \eqref{DE1E2} but without the factor $(2r - 1)(2s - 1)$ that appears in the second-order Peirce polynomial of $E_1 + 2E_2$. This makes it possible to calculate fusion rules for $O$ in the case when at least one of the eigenvalues is $\tfrac{1}{2}$. This means that the identity $E_1 + 2E_2$ gives no constraints on fusion rules involving $\frac{1}{2}$ even though the identity $O$ does. This difference in fusion rules is typical for the classes of algebras satisfying a univariate identity $P$ (such as $E_1 + 2E_2$) and the class of algebras satisfying the multivariate identity given by some linearization of $P$ (such as $O$). The latter class is contained in the former and is subject to more stringent fusion rules. \qed
\end{remark}

Just as the vanishing of the first order Peirce polynomial of an identity $P$ means that the abstract Peirce spectrum of $P$ is arbitrary, the vanishing of the second order Peirce polynomial of an evanescent identity $P$ means that the abstract fusion rules of $P$ are arbitrary.

\begin{definition}\label{def:secondorderevanescence}
A nonassociative polynomial $P\in \Nring{\field}{\{x\}}$ is \emph{second-order evanescent} if $P$ is evanescent and $\frak{D}(P;q, r, s)$ is identically zero.
\end{definition}

\begin{lemma}
The subspace $\evan^{(2)} := \ker \rho \cap \ker \mathfrak{D} \subset \Nring{\field}{\{x\}}$ of second-order evanescent polynomials is an ideal.
\end{lemma}

\begin{proof}
This follows from \eqref{secondpeircerecursion} of Lemma \ref{secondpeircerecursionlemma} and that, by Lemma \ref{evanescentp1lemma}, for $P \in \evan$ there holds $P(1) = 0$.
\end{proof}

By Proposition \ref{pro:de1e2}, the ideal of second order evanescent polynomials $\evan^{(2)}$ is a proper ideal in $\evan$. Example \ref{ex:secondorderevanescence} shows that $\evan^{(2)}$ is not the zero ideal.

\begin{example}\label{ex:secondorderevanescence}
Using Corollary~\ref{cor:Devan}, it follows that the nontrivial nonassociative polynomial
\begin{equation}\label{Pi}
\begin{aligned}
G&:=[E_1, x, x^2, x^2] \\
&= \big((x^2)^3x\big)(x^2)^2 -
\big((x^3)^2 x\big)(x^2)^2 - (x^2)^4x^3 + ((x^3)^2x^2)x^3
\end{aligned}
\end{equation}
satisfies
\begin{align}\label{Dzero}
&\rho(G;\lambda)=0,&
&\frak{D}(G; \alpha, \lambda, \mu)=0.
\end{align}
Let $\evan^{(2)}_{k}\subset \evan^{(2)}$ denote the homogeneous degree $k$ second-order evanescent polynomials.
This $G$ gives an example of a nonzero element of $\evan^{(2)}_{11}$. It follows from Proposition \ref{pro:de1e2} that no nonzero degree $6$ evanescent polynomial is second order evanescent, $\evan^{(2)}_{6} = \{0\}$.
It is not known what is the minimal degree $k$ such that $\evan^{(2)}_{k}\subset \evan^{(2)}$ is nonzero, but it has been established that such $k$ satisfies $7 \leq k \leq 11$.
\qed\end{example}

For particular choices of $(t_1, t_2)$, the discriminant $\Delta_{(t_1, t_2)}$ in $q$ of $\frak{D}(t_1 E_1 + t_2 E_2; q, \lambda, \mu)$,
\begin{align}
\begin{aligned}
\Delta_{(t_1, t_2)}
& = 4t_1^2 + 4\left(4\lambda\mu - \lambda - \mu - 2\right)t_1 t_2 + \left[8(\lambda - \mu)^2 + (\lambda + \mu - 2)^2\right]t_2^2,
\end{aligned}
\end{align}
has special features. This gives an alternative way of distinguishing identities within $\evan_6$. As an illustration, the fusion rules for certain of the corresponding evanescent identities are described briefly below.

\begin{example}[Fusion rules for $E_2$]
When $t_1 = 0$ and $t_2 = 1$, the resulting discriminant is a sum of three squares
\begin{equation}
\Delta_{(0,1)}
=4(\lambda - \mu)^2 + 4(\lambda - \mu)^2  + (\lambda + \mu - 2)^2.
\end{equation}
Straightforward calculations using \eqref{DE1E2} imply fusion rules for $E_2$ which are more complicated than those \eqref{peirceL} obtained later for $(2,3)$-palintropic algebras. If $\lambda$ and $\mu$ are both distinct from $\tfrac12$ then
\begin{equation}\label{Sprime0}
\alg_c(\lambda)\alg_c(\mu)\subset \alg_c(\alpha)\oplus\alg_c(\beta),
\end{equation}
provided $\alpha \neq \beta$ or $L(c)$ is diagonalizable over $\field$, 
where $\alpha$ and $\beta$ are the solutions of
\begin{align}\label{Sprimealpha1}
\begin{aligned}
&2(\alpha+\beta)=\lambda+\mu+2,&\\
&4(\alpha-\beta)^2=8(\lambda-\mu)^2+(\lambda+\mu-2)^2 = 9(\lambda - \mu)^{2} + 4(\lambda -1)(\mu-1),
\end{aligned}
\end{align}
provided these equations have solutions in $\field$. In particular, if the ground field is $\fie = \mathbb{R}$, then \eqref{Sprimealpha1} always has real solutions $\alpha, \beta$.

If $\lambda = \mu$ then the solutions of \eqref{Sprimealpha1} are $\{\alpha, \beta\} = \{\lambda, 1\}$.
If moreover $\lambda \notin \{\tfrac{1}{2}, 1\}$, these solutions are distinct elements of $\field$, and the fusion rule \eqref{Sprime0} becomes
\begin{align}\label{Sprime2}
\alg_c(\lambda)\alg_c(\lambda) \subset \alg_c(1) \oplus \alg_c(\lambda), \qquad \lambda \notin \{\tfrac{1}{2}, 1\}.
\end{align}
For $\lambda = \mu = 1$ the two solutions coincide, $\alpha = \beta = 1$, and $\frak{D}(E_2; q, 1, 1) = 2(q - 1)^2$, so \eqref{fusionrules0} yields only
\begin{align}\label{E2oneonefusion}
\alg_c(1)\alg_c(1) \subset \ker\bigl(L(c) - \Id\bigr)^2.
\end{align}
In particular, $\alg_c(1)$ is a subalgebra whenever $\ker(L(c) - \Id)^2 = \alg_c(1)$, for example when $L(c)$ is diagonalizable.
Similarly, if $\mu = 1$ then the solutions of \eqref{Sprimealpha1} are $\{\alpha, \beta\} = \{\lambda, \tfrac{1}{2}(3 - \lambda)\}$, which are distinct if and only if $\lambda \neq 1$. If moreover $\lambda \notin \{\tfrac{1}{2}, 1\}$, the fusion rule \eqref{Sprime0} becomes
\begin{align}\label{Sprime3}
\alg_c(\lambda)\alg_c(1) \subset \alg_c(\lambda) \oplus \alg_c(\tfrac{1}{2}(3 - \lambda)), \qquad \lambda \notin \{\tfrac{1}{2}, 1\}.
\end{align}
Although the abstract Peirce spectrum $\sigma_{\field}(E_2) = \field$ is unconstrained, the abstract fusion rules $\fusion_{\field}(E_2)$ are not: if $\lambda, \mu \neq \tfrac{1}{2}$ are eigenvalues of $L(c)$ and the product $\alg_c(\lambda)\alg_c(\mu)$ is nonzero, then at least one of the solutions $\alpha, \beta$ of \eqref{Sprimealpha1} lies in $\field$ and is an eigenvalue of $L(c)$. This implicitly puts constraints on the interactions between Peirce eigenvalues that appear to be quite subtle and which are not analyzed here. A similar remark applies to other identities $E$ in $\evan_{6}$.
\qed
\end{example}

\begin{example}[$E_4:=2E_1+E_2$]
The case $(t_1, t_2) = (2, 1)$ corresponding with the evanescent identity
\begin{align*}
E_4&:=2E_1+E_2=2 (x^2)^3 - x^3x^3 + x(x(x^2)^2)-x^2x^4 - x(x^2x^3),
\end{align*}
is distinguished because in this case the discriminant is a square of a linear form:
\begin{equation}
\Delta_{(2,1)} = (3\lambda + 3\mu - 2)^2.
\end{equation}
In this case,
$$
\frak{D}(E_4; q,\lambda, \mu) =   (2\lambda - 1)(2\mu - 1)(2q + \lambda + \mu)(q - \lambda - \mu + 1)
$$
and for $\lambda$ and $\mu$ distinct from $\tfrac12$ such that $\lambda + \mu \neq \tfrac23$ there hold the fusion rules
\begin{equation}\label{Sbis1}
\alg_c(\lambda)\alg_c(\mu)\subset \alg_c(-\tfrac12(\lambda + \mu))\oplus\alg_c(\lambda + \mu-1).
\end{equation}
If $\lambda + \mu = \tfrac23$, the two roots coincide, $\frak{D}(E_4;q,\lambda,\mu) = 2(2\lambda-1)(2\mu-1)(q+\tfrac13)^2$, and \eqref{fusionrules0} yields only $\alg_c(\lambda)\alg_c(\mu)\subset \ker(L(c)+\tfrac13\Id)^2$, which equals $\alg_c(-\tfrac13)$ when $L(c)$ is diagonalizable over $\field$.
\qed
\end{example}

\section{$(2,3)$-palintropic algebras and their fusion rules}\label{sec:adiabatic}
This section explains a special feature of the fusion rules for $E_1$, which is that they impose a multiplicative structure on the Peirce subspaces.

\begin{theorem}
\label{pro:23algebra}
Suppose $\operatorname{char} \field \neq 2$ and $|\field| > 6$, and let $\alg$ be a $(2, 3)$-palintropic $\field$-algebra. For $c\in \Idm(\alg)$, if $\lambda, \mu\in \sigma(c)\setminus\{\tfrac12\}$, then
\begin{equation}\label{alphabeta}
\alg_c(\lambda)\alg_c(\mu)\subset \alg_c(\lambda\mu).
\end{equation}
Furthermore:
\begin{enumerate}
[label=$\mathrm{(S\arabic*)}$,itemsep=0.5ex,leftmargin=1cm]
\item\label{s:01}
$\alg_c(0)$ and $\alg_c(1)$ are subalgebras of $\alg$.
\item\label{s:02}
If $c$ is nonsingular, then $\alg_c(0) = \ker L(c)$ is an ideal of $\alg$.
\item\label{s:03} If $x \in \alg_c(\lambda)$ and $\lambda^{j} \neq \tfrac{1}{2}$ for $1 \leq j < \deg(x^{\delta})$, then
\begin{align}\label{s:03b}
c\,x^{\delta} = \lambda^{\deg(x^{\delta})}x^{\delta}.
\end{align}
In particular, if $\lambda^{j} \neq \tfrac{1}{2}$ for all $j \geq 1$, then \eqref{s:03b} holds for every nonassociative power $x^{\delta}$.
\end{enumerate}
\end{theorem}

\begin{proof}
Because $\operatorname{char} \field \neq 2$, the second linearization of \eqref{strange} yields
\begin{equation}\label{2ndlinear}
\begin{aligned}
0 & = (x^2z)(x^2y)+2(x^2y)(x(xz))+2(x^2z)(x(xy))+2x^3\bigl[y(zx)+z(yx)+x(yz)\bigr]\\
&\quad +4(x(xz))(x(xy))-4(x^2(xy))(xz)-4(x^2(xz))(xy)-4x^2((xy)(xz))\\
&\quad -(x^2)^2(yz)-2x^2(x^2(yz)).
\end{aligned}
\end{equation}
Specializing $x=c\in \Idm(\alg)$ in \eqref{2ndlinear} yields
\begin{equation}\label{ccyz}
\begin{aligned}
0 & = (cy)(cz)-2(cy)(c(cz))-2(cz)(c(cy))+4(c(cy))(c(cz)) \\
&\quad -c(yz) +2c(y(zc))+2c(z(yc)) -4c((cy)(cz))\\
& = (\Phi L(c)y)(\Phi L(c)z) - L(c)((\Phi y)(\Phi z)),
\end{aligned}
\end{equation}
where $\Phi(c) = \Id - 2L(c)$ is as in \eqref{Phidefined}.
If $y\in \alg_c(\lambda)$, $z\in \alg_c(\mu)$, then \eqref{ccyz} yields
\begin{align*}
(\lambda\mu-2\lambda\mu(\lambda+\mu)+4\lambda^2\mu^2)yz+\bigl(2(\lambda+\mu)-4\lambda\mu -1\bigr)c(yz)=0,
\end{align*}
so that
\begin{align}\label{peirceL}
(1-2\lambda)(1-2\mu)(c(yz)-\lambda \mu yz)=0.
\end{align}
The fusion rule \eqref{alphabeta} follows from \eqref{peirceL}. Claim \ref{s:01} follows immediately from \eqref{alphabeta}.

If $y \in \ker L(c)$, then \eqref{ccyz} yields that
\begin{align}\label{ideal1}
0 = c(yz) - 2c(y(zc)) = c(y((\Id - 2L(c))z))
\end{align}
for all $z \in \alg$. This shows that $(\ker L(c))(\im (\Id - 2L(c))) \subset \ker L(c)$. If $\tfrac12 \notin \sigma(c)$, then $\Id - 2L(c)$ is invertible, so $\im(\Id - 2L(c)) = \alg$ and $\ker L(c) = \alg_{c}(0)$ is an ideal, showing \ref{s:02}.

Finally, let $x \in \alg_c(\lambda)$. The claim \eqref{s:03} is proved by induction on $\deg(x^{\delta})$, the case $\deg(x^{\delta}) = 1$ being the hypothesis. Write $x^{\delta} = x^{\beta}x^{\gamma}$ with $|\beta|, |\gamma| < |\delta|$. By the inductive hypothesis, $x^{\beta} \in \alg_c(\lambda^{|\beta|})$ and $x^{\gamma} \in \alg_c(\lambda^{|\gamma|})$, and by assumption $\lambda^{|\beta|} \neq \tfrac{1}{2}$ and $\lambda^{|\gamma|} \neq \tfrac{1}{2}$. Hence \eqref{alphabeta} yields $x^{\delta} \in \alg_c(\lambda^{|\beta| + |\gamma|}) = \alg_c(\lambda^{|\delta|})$.
\end{proof}

\begin{corollary}\label{cor:lambdamedial}
If $\alg$ is a $(2, 3)$-palintropic algebra, $c\in \Idm(\alg)$ and $\lambda\in \sigma(c)$ satisfies $\lambda\notin\{0,\tfrac12\}$ and $\lambda^2 \neq \tfrac12$, then
\begin{equation}\label{lam1}
\lambda (x^2)^2-x^4=0\quad \text{for all  } x\in \alg_c(\lambda).
\end{equation}
\end{corollary}
\begin{proof}
Let $x \in \alg_c(\lambda)$. Because $\lambda \neq \tfrac{1}{2}$, \eqref{alphabeta} gives $x^2 \in \alg_c(\lambda^2)$, and because $\lambda^2 \neq \tfrac{1}{2}$, it gives $(x^2)^2 \in \alg_c(\lambda^4)$. Substituting $y = z = c$ in \eqref{2ndlinear} and using these facts yields
\begin{align}\label{quasipowerproof}
2\lambda(2\lambda - 1)\bigl(\lambda(x^2)^2 - x^4\bigr) = 0,
\end{align}
and the claim follows because $\lambda \notin \{0, \tfrac{1}{2}\}$.
\end{proof}

\begin{remark}
A property like that of Corollary \ref{cor:lambdamedial} holds in any medial algebra; see Proposition 6.3 and Remark~6.4 in \cite{Krasnov-Tkachev-medial}. One can think of \eqref{lam1} as an explicit quantitative deviation from the power-associativity. Note, however, that this identity does not hold on the entire algebra, but only on the Peirce subspace associated with $\lambda$. \qed
\end{remark}

\begin{proposition}\label{pro:unital}
A \emph{unital} $(2,3)$-palintropic algebra over a field $\field$ such that $\chr \field \nmid 30$ is power-associative.
\end{proposition}

\begin{proof}
A direct consequence of the linearization of \eqref{strange} is that for any $x$ in any $(2, 3)$-palintropic algebra $\alg$ there holds
\begin{align}\label{weird1}
L(x^3)L(x^2)+2L(x^3)L(x)^2- L((x^2)^2)L(x)-2L(x^2)^2L(x)&=0.
\end{align}
Let $e \in \alg$ be the unit, so that $\alg=\alg_e(1)$, i.e. $ex=x$ for all $x$. Applying \eqref{weird1} to $e$ yields
\begin{equation}\label{descentphi1}
x^3x^2-(x^2)^2x=0.
\end{equation}
Note that the identity \eqref{descentphi1} has degree less than the original identity \eqref{strange}. Linearizing \eqref{descentphi1} yields
\begin{equation}\label{descent2}
x^2(x^2y+2x(xy))+2(xy)x^3-(x^2)^2y-4(x^2(xy))x=0,
\end{equation}
and taking $y=e$ yields $(x^2)^2-xx^3=0$. By a theorem of Albert \cite{Albert-powerassociative}, over a field of characteristic prime to $30$ a commutative algebra satisfying $(x^2)^2-xx^3=0$ is power-associative, so $\alg$ is power-associative.
\end{proof}

In general if a unital commutative algebra satisfies a homogeneous identity $P$, then it satisfies the linearization of $P$ also; see \cite[Proposition 19]{varro2020peirce} for a precise statement. Corollary \ref{cor:subpa} and Proposition \ref{pro:pseudounital} illustrate this for unital $(2, 3)$-palintropic algebras and their subclass comprising algebras satisfying \eqref{pseudoJ} are Jordan.

\begin{corollary}\label{cor:subpa}
In a $(2,3)$-palintropic algebra $\alg$ over a field $\field$ such that $\chr \field \nmid 30$ and having nonzero idempotent $c$, the subalgebra $\alg_c(1)$ is power-associative.
\end{corollary}
\begin{proof}
Combine Proposition~\ref{pro:unital} with \ref{s:01} of Theorem~\ref{pro:23algebra}.
\end{proof}

\begin{proposition}\label{pro:pseudounital}
A unital algebra $\alg$ over a field $\field$ with characteristic prime to $6$ that satisfies \eqref{pseudoJ} is Jordan.
\end{proposition}

\begin{proof}
Linearizing \eqref{pseudoJ} in $x$ at $z$ and replacing $z$ by the unit $e$ yields the degree $5$ identity
\begin{equation}\label{lin02}
3x^2(x^2y)-2x^3(xy)-(x^2)^2y=0.
\end{equation}
Because $\operatorname{char}\field$ is prime to $6$, linearizing \eqref{lin02} in $x$ at $z$ and replacing $z$ by the unit $e$ yields the degree $4$ identity
\begin{equation}\label{lin03}
x(x^2y)-x^3y=0,
\end{equation}
linearizing \eqref{lin03} in $x$ at $z$ and replacing $z$ by the unit $e$ yields the degree $3$ identity
\begin{equation}\label{lin04}
x(xy)-x^2y=0.
\end{equation}
Multiplying \eqref{lin04} by $x$ yields
\begin{equation}\label{lin05}
x(x(xy))-x(x^2y)=0.
\end{equation}
On the other hand, replacing $y$ in \eqref{lin04} by $xy$ yields
\begin{equation}\label{lin06}
x(x(xy))-x^2(xy )=0.
\end{equation}
Combining \eqref{lin05} and \eqref{lin06} implies $x^2(xy)=x(x^2y)$.
This shows the claim.
\end{proof}

\section{Almost mediality}
It is unreasonable to expect that the class of $(2,3)$-palintropic algebras is exhausted by the list in Proposition \ref{pro:classes},
even in the absence of the universal Peirce eigenvalue $\tfrac12$.
By the main theorem of \cite{Kokoris}, combined with results from \cite{Albert-powerassociative}, a simple finite-dimensional commutative power associative algebra over a field of characteristic zero that contains a nonzero idempotent is Jordan. However, in general Jordan algebras, defined by the two-variable identity~\eqref{jordandef}, are a proper subclass of power associative algebras, defined by the univariate identity~\eqref{pow}, \cite{rodrigues2020commutative}. Remark \ref{remark:Osborn} observed that something similar happens in regards to the $E_1 + 2E_2$ identity and algebras satisfying the Osborn identity.

Similarly, $(2, 3)$-palintropic algebras are more general than medial algebras.
Nevertheless, Theorem~\ref{the:iii} shows that certain properties of $(2,3)$-palintropic algebras are closely related to properties of medial algebras established in \cite{Tka23a}.

\begin{corollary}\label{almostmedialcorollary}
Suppose that $\alg$ is a $(2,3)$-palintropic algebra.
For $c\in \Idm(\alg)$ and $\Phi(c)= \Id - 2L(c)$, there holds
\begin{align}\label{prealmostmedial}
&(L(c)\Phi(c) y)(L(c)\Phi(c) z) = (\Phi(c) L(c)y)(\Phi(c) L(c)z) = L(c)((\Phi(c) y)(\Phi(c) z)), &
\end{align}
for all $y, z \in \alg$.
\end{corollary}

\begin{proof}
Because $\Phi$ and $L(c)$ commute, the last equality of \eqref{ccyz} yields \eqref{prealmostmedial}.
\end{proof}

It will be convenient to use \eqref{prealmostmedial} in the more explicit form
\begin{align}\label{almostmedial}
&(cy-2c(cy))(cz-2c(cz)) =c\left((y-2cy)(z-2cz)\right) & &\text{for all}\,\, y, z \in \alg.
\end{align}

\begin{corollary}\label{cor:iii}
Suppose $\chr \field \neq 2$.
Suppose that $\alg$ is a $(2, 3)$-palintropic $\field$-algebra. For $c\in \Idm(\alg)$ define $\Phi(c):=\Id_{\alg}-2L(c):\alg\to \alg$.
\begin{enumerate}
[label=$\mathrm{(\roman*)}$,itemsep=0.5ex,leftmargin=1cm]
\item For $\balg=\Phi(c)(\alg)$, $L(c):\balg\to \alg$ is a multiplicative linear map, meaning $L(c)(b_{1}b_{2}) = (L(c)b_{1})(L(c)b_{2})$ for $b_{1}, b_{2} \in \balg$.
\item Define $x \medmlt y := (\Phi(c)x)(\Phi(c)y)$. Then $L(c):(\alg, \medmlt) \to (\alg, \medmlt)$ is an algebra homomorphism.
\end{enumerate}
\end{corollary}

Note that in general it might be that $\balg$ is not a subalgebra.

\begin{proof}
Because $L(c)$ and $\Phi$ commute, after relabeling \eqref{ccyz} yields
\begin{align}
\begin{aligned}
L(c)(x\medmlt y) &= L(c)((\Phi(c)x)(\Phi(c)y))  = (L(c)\Phi(c)x)(L(c)\Phi(c)y) \\
&= (\Phi(c) L(c) x)(\Phi(c) L(c) y) = L(c)x \medmlt L(c)y,
\end{aligned}
\end{align}
from which both claims follow.
\end{proof}

Theorem \ref{the:iii} can be interpreted as showing that $(2, 3)$-palintropic algebras are almost medial because \cite[Proposition~6.1]{Tka23a} shows that any nonzero idempotent in a medial algebra satisfies the properties \ref{L:01}-\ref{L:03} of Theorem \ref{the:iii}. In the context of $(2, 3)$-palintropic algebras the same conclusion requires a nonsingular idempotent.

\begin{theorem}[Almost mediality]\label{the:iii}
If $c\in \Idm(\alg)$ is a nonsingular idempotent in the $(2,3)$-palintropic algebra $\alg$, then:
\begin{enumerate}
[label=$\mathrm{(M\arabic*)}$,itemsep=0.5ex,leftmargin=1cm]
\item\label{L:01}
$L(c):\alg \to \alg$ is an algebra homomorphism.
\item\label{L:02}
For any idempotent $c'\in\Idm(\alg)$, either $cc'=0$ or $cc'\in \Idm(\alg)$.
\item\label{L:03}
$\ker L(c)=\alg_c(0)$ is an algebra ideal.
\end{enumerate}
\end{theorem}

\begin{proof}
By assumption, $\Id -2L(c)$ is invertible. Given an arbitrary $z\in \alg$, let $y=(\Id -2L(c))^{-1}z$. Then $y-2cy=z$ and \eqref{almostmedial} yields
\begin{equation}\label{almostmedial2}
(cz)^2=cz^2, \qquad \text{for all}\,\, z\in \alg.
\end{equation}
Polarizing the identity \eqref{almostmedial2} yields
\begin{equation}\label{almostmedial3}
(cz)(cw)=c(zw), \qquad \text{for all}\,\, z,w\in \alg.
\end{equation}
In other words, $(L(c)z)(L(c)w)=L(c)(zw)$, so that $L(c)$ is an algebra homomorphism. This shows \ref{L:01}. Claim \ref{L:02} follows from \eqref{almostmedial2} and claim \ref{L:03} follows from \ref{L:01}.
\end{proof}

\begin{corollary}
A nonsingular idempotent in a finite-dimensional simple $(2, 3)$-palintropic algebra is invertible.
\end{corollary}
\begin{proof}
If $c$ is a nonsingular idempotent, then $c^{2} = c \neq 0$ so $\ker L(c) \neq \alg$. Because, by \ref{L:03} of Theorem \ref{the:iii}, $\ker L(c)$ is an ideal, this shows $L(c)$ is injective, and if the algebra is finite-dimensional this implies $c$ is invertible.
\end{proof}

\begin{remark}
Any Jordan algebra is $(2,3)$-palintropic, so Corollary~\ref{cor:iii} and Theorem \ref{the:iii} apply to Jordan algebras.
However, multiplication by an idempotent is not, in general, a Jordan algebra homomorphism.
To see what Corollary~\ref{cor:iii} implies in this case, recall that the Peirce decomposition of nonzero idempotent $c$ in a Jordan algebra $\alg$ has the form $\alg=\alg_c(1)\oplus\alg_c(0)\oplus\alg_c(\tfrac12)$. In the notation of
Corollary~\ref{cor:iii}, $\balg=\alg_c(0)\oplus \alg_c(1)$, where $\alg_c(0)$ and $\alg_c(1)$ are subalgebras of $\alg$, and
$$
L(c)(xy)=L(c)(x_0y_0+x_1y_1)=L(c)(x_1y_1)=x_1y_1=(L(c)x)(L(c)y)
$$ holds for $x, y \in \balg$.
However, $L(c)$ need not be a homomorphism of all of $\alg$: for $a, b \in \alg_c(\tfrac12)$ and $y \in \alg_c(0)$ there hold $c(ab) = (ab)_1$, $(ca)(cb) = \tfrac14 ab$, $c(ya) = \tfrac12 ya$ and $(cy)(ca) = 0$, so $L(c)$ fails to be a homomorphism as soon as $\alg_c(\tfrac12)\alg_c(\tfrac12) \neq 0$ (provided $\chr \field \neq 3$) or $\alg_c(0)\alg_c(\tfrac12) \neq 0$. For instance, in the Jordan algebra of symmetric $2\times 2$ matrices with product $x\circ y = \tfrac12(xy + yx)$, the idempotent $c = E_{11}$ and $a = E_{12} + E_{21} \in \alg_c(\tfrac12)$ satisfy
\begin{align}\label{jordannonhom}
c\circ(a\circ a) = E_{11} \neq \tfrac14 I = (c\circ a)\circ(c\circ a).
\end{align}
This shows that the assumption in Theorem \ref{the:iii} that the idempotent $c$ be nonsingular cannot be omitted.
\qed
\end{remark}

This section concludes with some examples of $(2, 3)$-palintropic algebras.

\smallskip
Suppose $p=(p_1,\ldots, p_n)\in \mathbb{Z}_+^n$, $p_i\ge2,$ and let $\field$ be a field of characteristic relatively prime to all the $p_i$.
Consider the commutative associative algebra $\balg=\field[u_1,\ldots,u_n]$ and its ideal $I_p:=\langle u_1^{p_1},\ldots,u_n^{p_n}\rangle$. Then the factor algebra $\balg/I_p$ is finite dimensional with $\dim \balg/I_p=p_1\ldots p_n$ and associative. Let $[\cdot ]:\balg\to \balg/I_p$ be the canonical projection. The linear endomorphism $\tau(u_1,\ldots, u_n)=(\lambda_1u_1,\ldots, \lambda_nu_n)$, where $\lambda_i\in \field$, stabilizes $I_p$, so descends to a linear endomorphism of the factor algebra $\balg/I_p$ defined by
$$
\tau([f(u)]):=[f(\tau(u))].
$$
 It follows from \cite{Tka23a} that the isotopic commutative algebra $(\alg,\ast)$ defined on $\alg=\balg/I_p$ by
\begin{equation}
\label{cutoff}
[f]\ast [g]=\tau([f(u)g(u)])=\tau([f(u)])\tau([g(u)])
\end{equation}
is a medial algebra, and thus is $(2,3)$-palintropic. Indeed, it follows from the definition that
\begin{equation}
\label{cutoff1}
([f]\ast [g])\ast ([h]\ast [k])=\tau^2([f(u)][g(u)][h(u)][k(u)])
\end{equation}
is completely symmetric.
This proves part of Proposition \ref{pro:severalvar}; the remaining claims are straightforward.

\begin{proposition}\label{pro:severalvar}
Given $p=(p_1,\ldots, p_n)\in \mathbb{Z}_+^n$,  a field $\field$ of characteristic coprime with any $p_i$ and $\tau(u_1,\ldots, u_n)=(\lambda_1u_1,\ldots, \lambda_nu_n)$, where $\lambda_i\in \field$,  let $(\alg,\ast)$  denote the isotope of the commutative algebra $\field[u_1,\ldots,u_n]/\langle u_1^{p_1},\ldots,u_n^{p_n}\rangle$ with the  multiplication \eqref{cutoff}. Then
\begin{enumerate}
[label=$\mathrm{(P\arabic*)}$,itemsep=0.5ex,leftmargin=1cm]
\item $(\alg,\ast)$ is a medial algebra of dimension $\dim \alg=p_1\ldots p_n$.
  \item $e=1\in \Idm(\alg,\ast)$.
  \item $u_i^{\ast k}=0$ for all $k\ge p_i$.
  \item $e\ast u_i^k=\lambda_i^k u_i^k$, i.e. $\field\{u_i^k\}\subset \alg_{e}(\lambda_i^k)$, $1\le k\le p_i-1$.
  \item In general, $\field\{u_1^{k_1}\cdots u_n^{k_n}\}\subset \alg_{e}(\lambda_1^{k_1}\cdots \lambda_n^{k_n})$, $0\le k_i\le p_i-1$.
\end{enumerate}
\end{proposition}

The construction in Proposition~\ref{pro:severalvar} is a variant of that used in \cite[Theorem~$8.8$]{Tka23a} to construct a diversity of medial algebras which relies on a different dilation action on the ideal $I_p$ and which necessarily entails the existence of nilpotent elements. The variation has the consequence that any finite set of elements of $\field$ can be realized in the Peirce spectrum of an idempotent in an algebra as in Proposition~\ref{pro:severalvar}.

\begin{example}
\label{ex:3}
This gives examples of $(2,3)$-palintropic algebras that are not medial and shows that the absence of the Peirce eigenvalue $\tfrac12$ is essential for the validity of the multiplicative property of Peirce eigenspaces \eqref{alphabeta}.

Let $\alg$ be the $\field$-algebra with basis $\{e_{i}: 0 \leq i \leq 3\}$, satisfying the relations
\begin{align*}
e_0^2&=e_0, && e_0e_i=\lambda_ie_i,\quad 1\le i\le 2,\\
e_0e_3&=p\lambda_1\lambda_2e_3, &&e_1e_2=qe_3,
\end{align*}
for $\lambda_1,\lambda_2, p, q\in \field$, with all other products of basis elements equal to $0$. Then for any $x=x_0e_0+x_1e_1+x_2e_2+x_3e_3$,
\begin{equation}
(x^2)^3-x^3x^3=2\lambda_1\lambda_2(2\lambda_1 - 1)(2\lambda_2 - 1)q(p - 1)x_0^4x_1x_2e_3
\end{equation}
It follows that $\alg$ is a $(2, 3)$-palintropic algebra whenever the condition
\begin{equation}\label{cond1}
2\lambda_1\lambda_2(2\lambda_1 - 1)(2\lambda_2 - 1)q(p - 1)=0
\end{equation}
is satisfied. In particular, consider the algebra satisfying
\begin{align}
\label{algnonmed}
\begin{aligned}
e_0^2&=e_0, && e_0e_1=\tfrac12 e_1,\\
e_0e_2&=\lambda e_2, && e_0e_3=\mu e_3,\qquad \mu\ne \tfrac12\lambda\\
e_1e_2&=e_3.
\end{aligned}
\end{align}
Because
$$
(e_0e_1)(e_0e_2)=\tfrac12\lambda e_1e_2=\tfrac12\lambda e_3\ne
(e_0e_0)(e_1e_2)=e_0e_3=\mu e_3,
$$
the algebra of \eqref{algnonmed} is not medial.

In this case, \eqref{alphabeta} fails for $\lambda_1=\tfrac12$ and $\lambda_2=\lambda$. Precisely,
$$
\alg_{e_0}(\tfrac12) \alg_{e_0}(\lambda)=\alg_{e_0}(\mu)\not\subset \alg_{e_0}(\tfrac12 \lambda)
$$
This shows that the restriction on the spectrum of $c$ in Theorem~\ref{pro:23algebra} is essential.
\qed\end{example}

\section{Third-order linearization of the $(2,3)$-palintropic identity}
The third order linearization of an identity yields information about threefold products of elements of an algebra satisfying the identity. Theorem \ref{pro:23trilinear} illustrates how the third order linearization yields information about $(2, 3)$-palintropic algebras. Theorem \ref{th:power} shows that in a $(2, 3)$-palintropic algebra, except for certain excluded eigenvalues $\lambda$, any monomial power of an element of the $\lambda$ Peirce eigenspace is proportional to a principal power.

\begin{theorem}
\label{pro:23trilinear}
Suppose $\operatorname{char}\field \notin \{2, 3\}$, let $\alg$ be a $(2,3)$-palintropic $\field$-algebra, and let $c \in \Idm(\alg)$. Let $x_i \in \alg_c(\lambda_i)$, $1 \leq i \leq 3$, where $\lambda_i \neq \tfrac{1}{2}$ and $\lambda_i\lambda_j \neq \tfrac{1}{2}$ for $1 \leq i < j \leq 3$. Then
\begin{equation}\label{eq:master_setting1}
W(x_1, x_2, x_3) :=w_1(\lambda_1, \lambda_2, \lambda_3)\xi_1 + w_2(\lambda_1, \lambda_2, \lambda_3)\xi_2 + w_3(\lambda_1, \lambda_2, \lambda_3)\xi_3 = 0,
\end{equation}
where $\xi_1 = (x_2 x_3)x_1,$ $\xi_2 = (x_3 x_1)x_2,$ $\xi_3 = (x_1 x_2)x_3$ and
\begin{align}\label{eq:masteri}
w_i = 4 \lambda_i \left( 2\lambda_i - \lambda_j - \lambda_k - 2\lambda_i\lambda_j - 2\lambda_i\lambda_k + 4\lambda_j\lambda_k \right), \qquad \{i,j,k\}=\{1,2,3\}.
\end{align}
\end{theorem}

\begin{proof}
Applying the multiplicative fusion rule \eqref{alphabeta} repeatedly yields
\begin{align}\label{rules}
\begin{aligned}
cx_k &= \lambda_k x_k \\
c(x_i x_j) &= \lambda_i \lambda_j (x_i x_j)\quad \text{for any  }  i\ne j, \quad i,j \in \{1,2,3\} \\
c((x_i x_j)x_k) &= \lambda_1 \lambda_2 \lambda_3 (x_i x_j)x_k \quad \text{for any permutation } (i,j,k) \text{ of } (1,2,3).
\end{aligned}
\end{align}
By \eqref{alphabeta}, $x_i x_j \in \alg_c(\lambda_i\lambda_j)$, and since $\lambda_i\lambda_j \neq \tfrac{1}{2}$ and $\lambda_k \neq \tfrac{1}{2}$, a second application of \eqref{alphabeta} gives $(x_i x_j)x_k \in \alg_c(\lambda_1\lambda_2\lambda_3)$.

Linearizing the identity \eqref{strange} defining $(2, 3)$-palintropic algebras three times yields $S_{1} - S_{2} = 0$ where
$S_1$ and $S_2$ are the third order linearizations of $(x^2)^3$ and $(x^3)^2$ given by
\begin{align}
\begin{aligned}
S_{1} &= 4\sum_{\sigma \in S_{3}}\left((x_{\sigma(1)}x_{\sigma(2)})(x_{\sigma(3)}y)\right)y^{2} + 4\sum_{\sigma \in S_{3}}\left((x_{\sigma(3)}y)y^{2}\right)(x_{\sigma(1)}x_{\sigma(2)})\\
&+ 4\sum_{\sigma \in S_{3}}\left( 2(x_{\sigma(1)}y)(x_{\sigma(2)}y) + (x_{\sigma(1)}x_{\sigma(2)})y^{2}\right)(x_{\sigma(3)}y), &\\
S_2 &= 2\sum_{\sigma \in S_{3}}\left((x_{\sigma(1)}x_{\sigma(2)})x_{\sigma(3)}\right)y^{3} \\
&+ 2\sum_{\sigma \in S_{3}}\left((x_{\sigma(1)}x_{\sigma(2)})y + (x_{\sigma(2)}y)x_{\sigma(1)} + (yx_{\sigma(1)})x_{\sigma(2)}\right)\left(x_{\sigma(3)}y^{2} + 2y(yx_{\sigma(3)})\right).
\end{aligned}
\end{align}
Taking $y = c \in \Idm(\alg)$ and evaluating $S_1$ and $S_2$ using \eqref{rules} yields
\begin{align*}
S_1 = \; &8\lambda_1\left(\lambda_1\lambda_2\lambda_3 + 3\lambda_2\lambda_3 + \lambda_1\right) (x_2x_3)x_1\\
    + \; &8\lambda_2\left(\lambda_1\lambda_2\lambda_3 + 3\lambda_1\lambda_3 + \lambda_2\right) (x_1x_3)x_2 \\
    + \; &8\lambda_3\left(\lambda_1\lambda_2\lambda_3 + 3\lambda_1\lambda_2 + \lambda_3\right) (x_1x_2)x_3
\end{align*}
and
\begin{align*}
S_2 = \; &4\lambda_1\left(2\lambda_1\lambda_2\lambda_3 + 2\lambda_1\lambda_2 + 2\lambda_1\lambda_3 + 2\lambda_2\lambda_3 + \lambda_2 + \lambda_3\right) (x_2x_3)x_1 \\
    + \; &4\lambda_2\left(2\lambda_1\lambda_2\lambda_3 + 2\lambda_1\lambda_2 + 2\lambda_1\lambda_3 + 2\lambda_2\lambda_3 + \lambda_1 + \lambda_3\right) (x_1x_3)x_2 \\
    + \; &4\lambda_3\left(2\lambda_1\lambda_2\lambda_3 + 2\lambda_1\lambda_2 + 2\lambda_1\lambda_3 + 2\lambda_2\lambda_3 + \lambda_1 + \lambda_2\right) (x_1x_2)x_3.
\end{align*}
Substituted in $0 = S_1 - S_2$, these yield \eqref{eq:master_setting1}.
\end{proof}

\begin{remark}
Note that by Theorem~\ref{pro:23algebra}, all $\xi_i$ lie in the same Peirce eigenspace $\alg_c(\lambda_1\lambda_2\lambda_3)$. Furthermore, if at least some of $\lambda_1, \lambda_2, \lambda_3$ are distinct then \eqref{eq:master_setting1} shows that the $\xi_i$ are linearly dependent. \qed
\end{remark}

The formulation of Lemma \ref{lem:power} excludes $\lambda$ contained in the set defined in \eqref{LAMBDA}. It is not clear whether this exclusion is essential for the validity of the main result below or merely reflects a limitation of the method of proof. Define $\Lambda$ as the set of all $\lambda \in \field$ for which there exist integers $k \geq 2$ and $n \geq 1$ such that
\begin{align}\label{LAMBDA}
R_{k,n}(\lambda) := 2\lambda^{k+n} - 4\lambda^{n+1} + 2\lambda^{k+1} + \lambda^{n} - 2\lambda^{k} + \lambda = 0.
\end{align}
For $k \geq 2$, $R_{k,n}$ is a nonzero polynomial: its term of highest degree is $2\lambda^{k+n}$ if $n \geq 2$ and $4\lambda^{k+1}$ if $n = 1$. Note that
\begin{align}\label{Rkn}
R_{k,k}(\lambda) = \lambda(\lambda^{k-1} - 1)(2\lambda^{k} - 1), \qquad R_{k,1}(\lambda) = 2\lambda(2\lambda - 1)(\lambda^{k-1} - 1).
\end{align}
Consequently $0$, $\tfrac{1}{2}$ and $1$ lie in $\Lambda$, and every $\lambda \notin \Lambda$ satisfies $\lambda^{m} \neq \tfrac{1}{2}$ for all $m \geq 1$ (take $R_{2,1}$ for $m = 1$ and $R_{m,m}$ for $m \geq 2$).

\begin{remark}\label{rem:Lambdaroots}
The set $\Lambda$ contains every root of unity in $\field$. Indeed, if $\lambda \in \field$ satisfies $\lambda^{m} = 1$ for some integer $m \geq 1$, then \eqref{Rkn} with $k = m+1 \geq 2$ gives
\begin{align}\label{Rrootsofunity}
R_{m+1,m+1}(\lambda) = \lambda(\lambda^{m} - 1)(2\lambda^{m+1} - 1) = 0,
\end{align}
so $\lambda \in \Lambda$. Since moreover $0 \in \Lambda$, and every nonzero element of a finite field is a root of unity, $\Lambda = \field$ whenever $\field$ is finite. Consequently, Lemma \ref{lem:power} and Theorem \ref{th:power} have content only when $\field$ is infinite. \qed
\end{remark}

\begin{lemma}
\label{lem:power}
Suppose $\operatorname{char} \field \notin \{2,3\}$, and let $\alg$ be a $(2, 3)$-palintropic $\field$-algebra.
Suppose $\lambda\not \in \Lambda$, $c\in \Idm(\alg)$ and $x\in \alg_c(\lambda)$. Then for any positive integers $k,n$,
\begin{equation}\label{knmain1}
x^nx^k=\lambda^{-(n-1)(k-1)}x^{n+k}.
\end{equation}
\end{lemma}

\begin{proof}
The lemma is proved by induction. Let $C(n)$ denote the statement that the conclusion of the lemma holds for $n$ and for all integers $k\ge1$. For $n=1$ the statement $C(1)$ asserts that $x^1x^k= x^{1+k}$, which holds by definition of principal powers.

For $n \geq 1$, suppose that $C(i)$ holds for $1\le i\le n$. It must be shown that $x^{n+1}x^{k}=\lambda^{-n(k-1)}x^{n+k+1}$ for all $k\geq 1$. For $k=1$ this reads $x^{n+1}x=x^{n+2}$, which holds by definition of principal powers. Let $k\geq 2$. Because $\lambda\notin\Lambda$, \eqref{Rkn} and $R_{2,1}(\lambda)=2\lambda(2\lambda-1)(\lambda-1)$ show that $\lambda^{m}\neq\tfrac12$ for all $m\geq1$. Hence \ref{s:03} of Theorem \ref{pro:23algebra} gives $x^{n}\in\alg_c(\lambda^{n})$ and $x^{k}\in\alg_c(\lambda^{k})$, and the hypotheses of Theorem \ref{pro:23trilinear} are satisfied for $x_1=x$, $x_2=x^n$, $x_3=x^k$ and $\lambda_1=\lambda$, $\lambda_2=\lambda^n$, $\lambda_3=\lambda^k$, since $\lambda$, $\lambda^{n}$, $\lambda^{k}$, $\lambda^{n+1}$, $\lambda^{k+1}$ and $\lambda^{n+k}$ all differ from $\tfrac12$. Then \eqref{eq:master_setting1} yields
\begin{align}\label{knmain2}
w_1(\lambda,\lambda^n,\lambda^k)\, x(x^nx^k)+w_2(\lambda,\lambda^n,\lambda^k)\, x^n(x^kx)+w_3(\lambda,\lambda^n,\lambda^k)\, x^{n+1}x^k=0.
\end{align}
By the inductive assumption $C(n)$,
\begin{align}\label{knmain3}
\begin{aligned}
x(x^nx^k)&=\lambda^{-(n-1)(k-1)}xx^{n+k}=\lambda^{-(n-1)(k-1)}x^{n+k+1},\\
x^n(x^kx)&=x^nx^{k+1}=\lambda^{-(n-1)k}x^{n+k+1},
\end{aligned}
\end{align}
which combined with \eqref{knmain2} yields
\begin{align}\label{knmain4}
-w_3(\lambda,\lambda^n,\lambda^k)\, x^{n+1}x^k=\bigl(w_1(\lambda,\lambda^n,\lambda^k) \lambda^{n-1} +w_2(\lambda,\lambda^n,\lambda^k)\bigr) \lambda^{-(n-1)k} x^{n+k+1}.
\end{align}
On the other hand, substituting $\lambda_1=\lambda$, $\lambda_2=\lambda^n$, $\lambda_3=\lambda^k$ in \eqref{eq:masteri} yields
\begin{align}\label{knmain5}
\begin{aligned}
w_1(\lambda,\lambda^n,\lambda^k) \lambda^{n-1} +w_2(\lambda,\lambda^n,\lambda^k)&=4\lambda^{n}R_{k,n}(\lambda),\\
-w_3(\lambda,\lambda^n,\lambda^k)&=4\lambda^{k}R_{k,n}(\lambda).
\end{aligned}
\end{align}
Because $\lambda\notin\Lambda$ and $k\geq 2$, $R_{k,n}(\lambda)\neq0$ by \eqref{LAMBDA}, and $\lambda\neq0$ because $0\in\Lambda$. Substituting \eqref{knmain5} in \eqref{knmain4} and dividing by $4\lambda^{k}R_{k,n}(\lambda)$ yields
\begin{align}\label{knmain6}
x^{n+1}x^{k}=\lambda^{n-k-(n-1)k}x^{n+k+1}=\lambda^{-n(k-1)}x^{n+k+1},
\end{align}
which is $C(n+1)$. This completes the induction step, and the lemma follows.
\end{proof}

It is sometimes convenient to write as $|\alpha|$ the degree $\deg(\alpha)$ of a nonassociative power $\alpha$, which is the number of leaves in the associated abstract binary rooted tree.
\begin{theorem}
\label{th:power}
Suppose $\operatorname{char} \field \notin \{2,3\}$, and let $\alg$ be a $(2, 3)$-palintropic $\field$-algebra.
For $\lambda\not \in \Lambda$, $c\in \Idm(\alg)$, and $x\in \alg_c(\lambda)$,
\begin{equation}\label{knmain}
x^\alpha =\lambda^{-\nu(x^\alpha)}x^{|\alpha|},
\end{equation}
where $|\alpha| = \deg(x^\alpha)$ and the map $\nu:\nass{\{x\}} \to \mathbb{Z}$ is uniquely determined by $\nu(x)=0$ and the recursion relation
\begin{equation}\label{knmain7}
\nu(x^\beta x^\gamma) =\nu(x^\beta)+\nu (x^\gamma)+(|\beta|-1)(|\gamma|-1).
\end{equation}
\end{theorem}

\begin{proof}
This follows from an induction based on
\begin{align}
\begin{aligned}
x^\beta x^\gamma&=\lambda^{-\nu(x^\beta)}x^{|\beta|}\cdot \lambda^{-\nu(x^\gamma)}x^{|\gamma|}=\lambda^{-\nu(x^\beta)-\nu(x^\gamma)}x^{|\beta|}x^{|\gamma|} \\
&=\lambda^{-\nu(x^\beta)-\nu(x^\gamma)}\cdot  \lambda^{-(|\beta|-1)(|\gamma|-1)}x^{|\beta|+|\gamma|}\\
&=\lambda^{-\nu(x^\beta x^\gamma)}x^{|\beta|+|\gamma|},
\end{aligned}
\end{align}
in which the third equality follows from Lemma~\ref{lem:power}.
\end{proof}
\begin{remark}
It is known \cite[Proposition~6.5]{Krasnov-Tkachev-medial} that for any medial algebra admitting a nonsingular invertible idempotent every monomial power $x^\alpha$ is proportional to the principal power $x^{\deg(\alpha)}$.
\end{remark}

Lemma \ref{lem:nuformula} gives a closed formula for $\nu$ in terms of $\mathfrak{D}$.

\begin{lemma}\label{lem:nuformula}
The map $\nu:\nass{\{x\}} \to \mathbb{Z}$ satisfies
\begin{equation} \label{H_nu_identity}
\nu(x^{\alpha}) = \tfrac{1}{2}(\deg(x^\alpha)+2)(\deg(x^\alpha) - 1) - \frak{D}(x^\alpha; 1, 1, \tfrac{1}{2}),
\end{equation}
where
\begin{equation} \label{H_nudef}
\frak{D}(x^\alpha; 1, 1, \tfrac{1}{2})
=\rho'(x^\alpha; 1),
\end{equation}
and $\rho'(x^\alpha; 1):=\tfrac{d}{dq}\rho(x^\alpha; q)|_{q=1}$.
\end{lemma}

\begin{proof}
Setting $H(x^{\alpha}):=\frak{D}(x^\alpha; 1, 1, \tfrac{1}{2})$, the goal is to prove that
\begin{equation} \label{H_nu_identity1}
H(x^{\alpha})=\tfrac12(\deg(x^\alpha)+2)(\deg(x^\alpha) - 1) -\nu(x^{\alpha}).
\end{equation}
The argument is by induction on the degree of monomials. For $x^\alpha = x$ the statement holds because
$$
\tfrac12(\deg(x^\alpha)+2)(\deg(x^\alpha) - 1) - \nu(x) = \tfrac12(1+2)(1-1) - 0 = 0 = H(x).
$$
Write $x^\alpha = x^\beta x^\gamma$ and let $n = \deg(x^\alpha)$, $k = \deg(x^\beta)$, and $m = \deg(x^\gamma)$, so $n = k + m$.
Expanding the components of $H(x^\alpha) = \frak{D}(x^\alpha; 1, 1, \frac{1}{2})$ using \eqref{defDD}(ii) yields the recursion:
\begin{equation}
\label{H_rec}
H(x^\beta x^\gamma) = H(x^\beta) + H(x^\gamma) + n = H(x^\beta) + H(x^\gamma) + k + m.
\end{equation}

The inductive hypothesis means the identity \eqref{H_nu_identity1} holds for $x^\beta$ and $x^\gamma$.
Substituting $H(x^\beta) = \tfrac12(k^2 + k - 2) - \nu(x^\beta)$ and $H(x^\gamma) = \tfrac12(m^2 + m - 2)-\nu(x^\gamma)$ into \eqref{H_rec} gives:
\begin{align*}
H(x^\beta x^\gamma)
&= \left( \tfrac12(k^2 + k - 2) - \nu(x^\beta)\right) + \left( \tfrac12(m^2 + m - 2) - \nu(x^\gamma)\right) + k + m \\
&= \tfrac12(k+m)^2 +  \tfrac12(k+m) - 1 - \nu(x^\beta) - \nu(x^\gamma) - (km - k - m + 1)\\
&= \tfrac12(n^2 + n - 2) - \left[ \nu(x^\beta) + \nu(x^\gamma) + (k-1)(m-1) \right]\quad \text{(by \eqref{knmain7})}\\
&= \tfrac12(n+2)(n-1) - \nu(x^\beta x^\gamma).
\end{align*}
This completes the induction and proves the main claim.

Formally, the relation \eqref{H_nudef} follows from setting $q$ and $r$ equal to $1$ in \eqref{Devan}; it can be justified by an induction starting from the recursion for $\frak{D}(x^{\alpha}, 1, 1, \tfrac{1}{2})$ obtained by setting $q = r = 1$ and $s = \tfrac{1}{2}$ in \eqref{defDD}.
\end{proof}

\begin{remark}
For the monomial $x^{\alpha} \in \Nring{\field}{\{x\}}$ associated with the abstract binary rooted tree $T$ having set of leaves $L(T)$, there holds the following explicit formula:
\begin{align}\label{derivativerhotreeformula}
\mathfrak{D}(x^{\alpha };1,1,\tfrac{1}{2})=\rho ^{\prime }(x^{\alpha };1)
=\sum _{v\in T\setminus L(T)}\deg(T_v),
\end{align}
where $T_v \subset T$ is the binary rooted subtree of $T$ with root $v$. The second equality in \eqref{derivativerhotreeformula} follows from differentiating at $q = 1$ the tree formula \eqref{peircetreeformula} for the Peirce operator of a monomial. 
The sum runs over the internal vertices of $T$; each leaf at height $h$ is counted once for each of its $h$ proper ancestors, which gives $\rho'(x^\alpha;1)=\sum_{v\in L(T)}\operatorname{ht}(v)$.
\end{remark}

\begin{example}
For the degree $5$ monomial $x(x^2x^2)$ shown in the diagram, $\rho(x(x^2x^2);q)=4q^3+q$, hence $\rho'(x(x^2x^2);1)=(12q^2+1)|_{q=1}=13$. On the other hand, calculating the degrees of the two base $x^2$ nodes ($\deg(T_v) = 2$), the intermediate $(x^2x^2)$ node ($\deg(T_v) = 4$), and the highest root node ($\deg(T_v) = 5$) yields:
\begin{equation*}
\frak{D}(x(x^2x^2);1, 1, \tfrac{1}{2}) = 5 + 4 + 2 + 2 = 13.
\end{equation*}
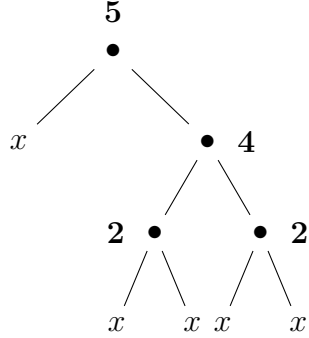
\begin{center}
\begin{tikzpicture}[
  level distance=1.2cm,
  level 1/.style={sibling distance=25mm},
  level 2/.style={sibling distance=14mm},
  level 3/.style={sibling distance=10mm}
]
  \node[ label=above:{\textbf{5}}] (v4) {$\bullet$}
    child {node {$x$}}
    child {node[ label=right:{\textbf{4}}] {$\bullet$}
      child {node[ label=left:{\textbf{2}}] {$\bullet$}
        child {node {$x$}}
        child {node {$x$}}
      }
      child {node[ label=right:{\textbf{2}}] {$\bullet$}
        child {node {$x$}}
        child {node {$x$}}
      }
    };
\end{tikzpicture}
\end{center}
In this case, $\nu(x(x^2x^2))=\tfrac12(5+2)(5-1)-13=1$  and \eqref{knmain} gives $x(x^2x^2)=\lambda^{-1}x^5$.
\qed
\end{example}

\section{Application to commuting and integrable polynomial maps}\label{sec:integrable}

The main point of this section is that commutative algebras satisfying evanescent identities of a certain kind yield examples of integrable polynomial mappings in the sense defined in \cite{veselov1991integrable}.

\begin{definition}\label{def:integrable}
A polynomial map of degree at least one $p:\field^n\to \field^n$ is \textit{integrable} if there exists a polynomial map $q:\field^n\to \field^n$ of degree at least one such that:
\begin{itemize}
\item $p\circ q=q\circ p$;
\item the sets of iterates of $p$ and $q$ are disjoint.
\end{itemize}
\end{definition}

Integrable maps play an important role in the theory of dynamical systems, as they exhibit an unusual degree of symmetry. For $n=1$, integrable polynomials were completely classified by Julia, Fatou, and Ritt in the early 1920s: any such map is linearly conjugate either to $z^n$ or to a Chebyshev polynomial; see \cite{Eremenko} for this and references. There is an analogous construction in higher dimensions,  related to Lie algebras via an action of an affine Weyl group \cite{kuccuksakalli2016bivariate, ChalykhVeselov}; for further information and references, see also \cite{Silverman}.

Given a commutative $\field$-algebra and a basis $\basis = \{e_{i}\}_{1 \leq i \leq n}$ of $\alg$, each $P \in \Nring{\field}{\{x\}}$ determines a polynomial map $\Phi(P):\field^{n}\to \field^{n}$ whose components $\Phi(P)(t)^{i}$ in the basis $\basis$ are given by
\begin{align}
\sum_{i = 1}^{n}\Phi(P)(t)^{i}e_{i} = P(\sum_{i = 1}^{n}t^{i}e_{i})
\end{align}
for $(t^{1}, \dots, t^{n}) \in \field^{n}$. Although $\Phi(P)$ depends on the choice of basis $\basis$, its equivalence class under conjugation by an element of $GL(n, \field)$ does not, and it is generally convenient to omit notation indicating dependence on $\basis$.

By Proposition \ref{pro:comp}, for any $P, Q \in \Nring{\field}{\{x\}}$ with $P(1)=Q(1)=1$, the identity $P\circ Q - Q\circ P$ is evanescent. If $\alg$ is a commutative $\field$-algebra satisfying $P\circ Q - Q\circ P$ then polynomial maps $\Phi(P)$ and $\Phi(Q)$ of $\field^{n}$ induced by $P$ and $Q$ and the choice of basis $\basis$ commute in the sense that
\begin{align}
\Phi(P)\circ \Phi(Q) = \Phi(Q)\circ \Phi(P).
\end{align}
However, it is not immediate that the sets of iterates of $\Phi(P)$ and $\Phi(Q)$ be disjoint, so it need not be the case that they are integrable polynomial mappings. On the other hand, for the simplest case where $P(x) = x^{\alpha}$ and $Q(x) = x^{\beta}$ are nonassociative monomials of relatively prime degrees, it is not hard to show that the resulting polynomials are integrable (provided no iterate vanishes) and this motivates focusing on that particular case in what follows.

This means that a commutative algebra satisfying $P\circ Q - Q\circ P$ and such that iterates of $\Phi(P)$ and $\Phi(Q)$ are disjoint can be viewed as a kind of \emph{quantization} of the integrable pair $\Phi(P)$ and $\Phi(Q)$. What is not obvious is what the disjoint iterates condition means at the level of $P$ and $Q$. Alternatively, an evanescent identity plus some condition generalizing the disjointness of the iterates of $\Phi(P)$ and $\Phi(Q)$ can be viewed as a kind of \emph{integrable nonassociative algebra}.

The remainder of this section focuses on the case of $(\alpha, \beta)$-palintropic algebras for abstract binary rooted trees $\alpha$ and $\beta$. In this case the nonassociative polynomials $P$ and $Q$ are the monomial power maps $P(x) = x^{\alpha}$ and $Q(x) = x^{\beta}$ and there are written $\Phi_{\alpha} = \Phi(P)$ and $\Phi_{\beta} = \Phi(Q)$ (for some choices of $\alg$ and basis $\basis$).  (In the case of the principal powers $P(x) = x^{p}$ and $Q(x) = x^{q}$, there are written $\Phi_{p} = \Phi(P)$ and $\Phi_{q} = \Phi(Q)$.) By assumption these satisfy $\Phi_\alpha \circ \Phi_\beta =\Phi_\beta \circ \Phi_\alpha$.
Proposition \ref{integrablepolynomialproposition} shows that these polynomial maps are integrable provided the degrees of $\alpha$ and $\beta$ are relatively prime.

Recall that $|\alpha|$ means $\deg(x^{\alpha})$.
\begin{proposition}\label{integrablepolynomialproposition}
Suppose $\field$ is infinite. Let $\alpha$ and $\beta$ be nonassociative powers of coprime degrees $|\alpha|, |\beta| \geq 2$, let $\alg$ be an $(\alpha, \beta)$-palintropic algebra with a fixed basis $\mathcal{B} = \{e_i\}_{1 \leq i \leq n}$, and let $\Phi_\alpha, \Phi_\beta : \field^n \to \field^n$ be the polynomial maps determined by $\alg$ and $\mathcal{B}$. Suppose that
\begin{align}\label{nonnilpotent}
\Phi_\alpha^{k} \not\equiv 0 \quad \text{and} \quad \Phi_\beta^{k} \not\equiv 0 \qquad \text{for all } k \geq 1.
\end{align}
Then $\Phi_\alpha$ and $\Phi_\beta$ are integrable mappings. Condition \eqref{nonnilpotent} holds whenever $\alg$ contains a nonzero idempotent.
\end{proposition}
\begin{proof}
Because $x \mapsto x^\alpha$ is homogeneous of degree $|\alpha|$, each component of $\Phi_\alpha^{k}$ is a homogeneous polynomial of degree $|\alpha|^{k}$, so by \eqref{nonnilpotent} $\Phi_\alpha^{k}$ is a nonzero homogeneous polynomial map of degree exactly $|\alpha|^{k}$; similarly $\Phi_\beta^{m}$ is a nonzero homogeneous polynomial map of degree exactly $|\beta|^{m}$. In particular $\Phi_\alpha$ and $\Phi_\beta$ have degree at least one. Since $\alg$ satisfies $(x^\alpha)^\beta = (x^\beta)^\alpha$ and $\field$ is infinite, $\Phi_\alpha \circ \Phi_\beta = \Phi_\beta \circ \Phi_\alpha$ as polynomial maps. Were $\Phi_\alpha^{k} = \Phi_\beta^{m}$ for some $k, m \geq 1$, comparing degrees would give $|\alpha|^{k} = |\beta|^{m}$, which is impossible because a prime dividing $|\alpha|$ would divide $|\beta|$. Hence the sets of iterates of $\Phi_\alpha$ and $\Phi_\beta$ are disjoint, and each map is integrable, the other serving as the commuting map in Definition \ref{def:integrable}.

Finally, if $c \in \Idm(\alg)$, an induction on the degree shows $c^\gamma = c$ for every nonassociative power $\gamma$, so $\Phi_\alpha^{k}(c) = c \neq 0$ and $\Phi_\beta^{k}(c) = c \neq 0$ for all $k \geq 1$, which gives \eqref{nonnilpotent}.
\end{proof}

\begin{example}\label{ex:2b}
Consider Proposition \ref{pro:severalvar} in the case $n=2$ and $p=(2,2)$ and let $\alg=(\field[u_1,u_2]/\langle u_1^2,u_2^2\rangle,\ast)$. Then $\dim \alg=2^2=4$. Consider the basis of $\alg$ given by $e_0=[1]$, $e_1=[u_1]$, $e_2=[u_2]$ and $e_3=[qu_1u_2]$, where $q = \lambda_1\lambda_2\in \dot{\field}$.
Then $\alg$ is the span of $e_{0}, e_{1}, e_{2}, e_{3}$ satisfying the relations
\begin{align}
\label{system1}
\begin{aligned}
&e_0^2=e_0, && \quad e_0e_i=\lambda_ie_i,\quad 1\le i\le 2,\\
&e_0e_3=\lambda_1\lambda_2e_3, &&\quad e_1e_2=e_3,
\end{aligned}
\end{align}
(all other products of basis elements are zero) where $\lambda_1,\lambda_2\in \dot{\field}$. Writing $x = \sum_{i = 0}^{3}x_{i}e_{i}$,
\begin{align}\label{ex:2product}
\begin{aligned}
xy &= x_{0}y_{0}e_{0} + \lambda_{1}(x_{0}y_{1} + x_{1}y_{0})e_{1} + \lambda_{2}(x_{0}y_{2} + x_{2}y_{0})e_{2}\\
&\quad + (\lambda_{1}\lambda_{2}(x_{0}y_{3} + x_{3}y_{0}) + x_{1}y_{2} + x_{2}y_{1})e_{3}.
\end{aligned}
\end{align}
After reindexing, this yields the commuting polynomial endomorphisms of $\field^4$:
\begin{align*}
P(t)&=( t_1^2,\, 2 \lambda_1 t_1 t_2  ,\,  2  \lambda_2 t_1 t_3,\,  2 \lambda_1  \lambda_2  t_1 t_4 + 2 t_2 t_3),\\
Q(t)&=(t_1^3,\,  \lambda_1   (2 \lambda_1  + 1)t_1^2t_2,  \lambda_2 (2 \lambda_2  + 1)t_1^2 t_3 ,   \lambda_1 \lambda_2 (2\lambda_1\lambda_2+1)  t_1^2 t_4 + 2 (\lambda_1  \lambda_2   + \lambda_1 + \lambda_2)  t_1 t_2t_3 ),
\end{align*}
where $t=(t_1,t_2,t_3,t_4)\in \field^4$.
Setting $t_1=1$ and $s_{i}=t_{i+1}$ for $i=1,2,3$, yields two commuting polynomial endomorphisms of $\field^3$:
\begin{align*}
p(s)&=(2\lambda_1   s_1  ,\,  2\lambda_2   s_2,\,  2\lambda_1  \lambda_2    s_3 +  2s_1 s_2),\\
q(s)&=(\lambda_1   (2 \lambda_1  + 1)s_1,  \lambda_2 (2 \lambda_2  + 1) s_2 ,   \lambda_1 \lambda_2 (2\lambda_1\lambda_2+1) s_3 + 2 (\lambda_1  \lambda_2   + \lambda_1 + \lambda_2)    s_1s_2 )
\end{align*}
where $s=(s_1,s_2,s_3)\in \field^3$.
Note that the latter polynomial maps are no longer homogeneous.
\qed\end{example}

A homogeneous quadratic polynomial map determines a commutative algebra whose idempotents are the fixed points of the map. Polarization associates with a homogeneous quadratic polynomial mapping $P:\field^n \to \field^n$ a commutative multiplication on $\field^n$ defined by
\begin{align}\label{polarquadratic}
x\star y=P(x,y):=\tfrac12(P(x+y)-P(x)-P(y)).
\end{align}
The relevant feature of the algebra $(\field^n, \star)$ is that, because $c \star c = P(c, c) = P(c)$, its idempotents are exactly the fixed points of $P$.

\begin{example}
In the context of commuting polynomial maps this yields the following interesting observation.
Suppose $P$ and $Q$ are commuting polynomial maps of homogeneous degrees $2$ and $q\ge2$, respectively, and define a commutative algebra $(\balg,\star)$ by polarizing $P$ as in \eqref{polarquadratic} above.
In this case, if $c$ is an idempotent in $(\balg, \star)$, then $Q(c)$ is a $\star$-idempotent too. Indeed, if $c = c\star c=P(c,c)=P(c)$, then $Q(c)=Q(P(c))=P(Q(c)) = Q(c)\star Q(c)$.
\qed\end{example}

Recall that a foldable  convex polytope in $\mathbb{R}^n$ can be subdivided into smaller congruent copies by reflections which corresponds bijectively to an affine Weyl group. For such foldable figure, one can construct a semigroup with a sequence $P_k(z)$ of polynomial maps of $\mathbb{R}^n$ satisfying the composition property
\begin{equation}\label{Pkm}
P_k\circ P_m=P_{km}.
\end{equation}
The Chebyshev polynomials $T_k(z)$ of the 1st kind, defined by $T_k(\cos\theta)=\cos k\theta$, constitute a well-known example of such sequences and, together with the trivial family $P_k(z)=z^k$ essentially exhaust all possibilities in dimension $n=1$.

In general, \eqref{Pkm} is extremely restrictive and characterizes these maps as higher-dimensional analogues of Chebyshev/power maps. Hoffman–Withers in \cite{HoffmanWithers} constructed, for every affine Weyl group, a family of polynomial maps $P_k$ that are conjugate to dilations and generalize Chebyshev polynomials via reflection symmetries, preserving composition, orthogonality, and geometric folding structure. In \cite{Withers} several explicit examples of such maps are given. The polynomial maps considered there exhibit properties of interest in the present nonassociative algebra context; they are described briefly now.

\begin{example}
Consider the sequence of maps $B_k$ defined by (7.1) in \cite[p.~409]{Withers}.
After a suitable homogenization, which increases the dimension by $1$, one can associate to the quadratic map $P_2$ a commutative algebra with multiplication defined by polarization as in \eqref{polarquadratic}. More precisely, $B_2(x,y)=(x^2-2y-4,y^2-2x^2+4y+4)$ can be homogenized as follows:
\begin{equation}\label{W2def}
W_2(x,y,z)=(x^2-2yz-4z^2,y^2-2x^2+4yz+4z^2,z^2):\field^3\to \field^3
\end{equation}
By \eqref{polarquadratic}, this gives the multiplication $\ast$ on $\field^3$ defined explicitly by
\begin{align}\label{B2mult}
\begin{aligned}
&(x_1,y_1,z_1)\ast (x_2,y_2,z_2) =\\
&=(x_1x_2-y_1z_2-y_2z_1-4z_1z_2,y_1y_2-2x_1x_2+2y_1z_2+2y_2z_1+4z_1z_2,z_1z_2)
\end{aligned}
\end{align}
The idempotents of the algebra $\balg:=(\field^3, \ast)$ are exactly the fixed points of \eqref{W2def}.
Let $(x,y,z)$ be a fixed point of \eqref{W2def}, then $z^2=z$. If $z=0$ then $x^2=x$ and $y^2-2x^2=y$, implying 4 fixed points
\begin{equation}\label{fixed1}
I_1=\{(0,0,0), \,(0,1,0), \,(1,-1,0), \,(1,2,0)\}.
\end{equation}
Similarly, if $z=1$ then  $x^2-2y-4-x=0$ and $y^2-2x^2+3y+4=0$, summing up obtain $(x-y)(x+y+1)=0$, thus implying 4 more fixed points
\begin{equation}\label{fixed2}
I_2=\{(1,-2,1), \,(-2,1,1), \,(-1,-1,1), \,(4,4,1)\}.
\end{equation}
The idempotents of $\balg$ are $\Idm(\balg)\cup \{0\} = I_1\cup I_2$. A further calculation yields the Peirce spectrum for each idempotent (in one case there must be supposed $\sqrt{-1} \in \field$):
\begin{align}
\begin{aligned}
&\sigma(0,0,0)=\{0,0,0\}, &
&\sigma(0,1,0)=\{0,0,1\},\\
&\sigma(1,-1,0)=\{-1,0,1\},&
&\sigma(1,2,0)=\{0,1,2\},\\
&\sigma(1,-2,1)=\{-1,1,2\},&
&\sigma(-1,-1,1)=\{1,\sqrt{-1},-\sqrt{-1}\},\\
&\sigma(-2,1,1)=\{-1,1,2\},&
&\sigma(4,4,1)=\{1,2,8\}.
\end{aligned}
\end{align}
Since there are exactly $2^3=8$ idempotents and $\tfrac12$ is absent in the Peirce spectrum, $\balg$ is a generic algebra according to the definition in \cite{Tkachev-universality}. This gives an example where 1) all idempotents can be calculated explicitly and they furthermore have integer coordinates and 2) the Peirce spectrum can be determined explicitly and consists only of (Gaussian) integers.
\qed\end{example}

\begin{small}
\def\polhk#1{\setbox0=\hbox{#1}{\ooalign{\hidewidth
  \lower1.5ex\hbox{`}\hidewidth\crcr\unhbox0}}} \def\cprime{$'$}
  \def\cprime{$'$} \def\cprime{$'$}
  \def\polhk#1{\setbox0=\hbox{#1}{\ooalign{\hidewidth
  \lower1.5ex\hbox{`}\hidewidth\crcr\unhbox0}}} \def\cprime{$'$}
  \def\cprime{$'$} \def\cprime{$'$} \def\cprime{$'$}
  \def\polhk#1{\setbox0=\hbox{#1}{\ooalign{\hidewidth
  \lower1.5ex\hbox{`}\hidewidth\crcr\unhbox0}}} \def\cprime{$'$}
  \def\Dbar{\leavevmode\lower.6ex\hbox to 0pt{\hskip-.23ex \accent"16\hss}D}
  \def\cprime{$'$} \def\cprime{$'$} \def\cprime{$'$} \def\cprime{$'$}
  \def\cprime{$'$} \def\cprime{$'$} \def\cprime{$'$} \def\cprime{$'$}
  \def\cprime{$'$} \def\cprime{$'$} \def\cprime{$'$} \def\dbar{\leavevmode\hbox
  to 0pt{\hskip.2ex \accent"16\hss}d} \def\cprime{$'$} \def\cprime{$'$}
  \def\cprime{$'$} \def\cprime{$'$} \def\cprime{$'$} \def\cprime{$'$}
  \def\cprime{$'$} \def\cprime{$'$} \def\cprime{$'$} \def\cprime{$'$}
  \def\cprime{$'$} \def\cprime{$'$} \def\cprime{$'$} \def\cprime{$'$}
  \def\cprime{$'$} \def\cprime{$'$} \def\cprime{$'$} \def\cprime{$'$}
  \def\cprime{$'$} \def\cprime{$'$} \def\cprime{$'$} \def\cprime{$'$}
  \def\cprime{$'$} \def\cprime{$'$} \def\cprime{$'$} \def\cprime{$'$}
  \def\cprime{$'$} \def\cprime{$'$} \def\cprime{$'$} \def\cprime{$'$}
  \def\cprime{$'$} \def\cprime{$'$}

\end{small}

\end{document}